\documentclass[11pt,reqno]{amsart}

\usepackage[utf8]{inputenc} 
\usepackage{comment}
\usepackage{amsmath}
\usepackage{amssymb}
\usepackage{amsfonts}
\usepackage{amsthm}
\usepackage{multirow}
\usepackage{enumerate}
\usepackage{color}
\usepackage{xcolor}
\usepackage{booktabs}
\usepackage{braket}
\usepackage[a4paper,margin=30mm]{geometry}
\usepackage{microtype}
\usepackage{hyperref}

\hypersetup{colorlinks=true,linkcolor=blue!40!black,citecolor=blue!40!black,urlcolor=blue!40!black}
\newcommand{\Z}{\mathbb{Z}}
\newcommand{\N}{\mathbb{N}}
\newcommand{\R}{\mathbb{R}}
\newcommand{\C}{\mathbb{C}}
\newcommand{\K}{K{\"a}hler }
\newcommand{\GL}{\operatorname{GL}}
\newcommand{\Id}{\operatorname{Id}}
\newcommand{\Aut}{\operatorname{Aut}}
\newcommand{\Vol}{\operatorname{Vol}}
\newcommand{\rank}{\operatorname{rank}}
\newcommand{\Ric}{\operatorname{Ric}}
\newcommand{\grad}{\operatorname{grad}}
\newcommand{\Trace}{\operatorname{Trace}}
\newcommand{\Jac}{\operatorname{Jac}}

\newcommand{\logterm}{{\rm{logterm}}}

\newcommand{\FS}{\mathrm{FS}}

\newtheorem{theor}{Theorem}[section]
\newtheorem{prop}[theor]{Proposition}

\newtheorem{cor}[theor]{Corollary}

\newtheorem{rmk}[theor]{Remark}

\title[The Szeg\H{o} log-term and local rigidity]{The Szeg\H{o} log-term and local rigidity of K\"ahler--Einstein polarizations of compact Hermitian symmetric spaces}

\author[A. Loi]{Andrea Loi}
\email{loi@unica.it}
\author[M. Palmieri]{Matteo Palmieri}
\email{matteo.palmieri@unica.it}
\author[F. Zuddas]{Fabio Zuddas}
\email{fabio.zuddas@unica.it}


\date{}

\begin{document}

\begin{abstract}
We prove local rigidity of the vanishing Szeg\H{o} log-term for integral K\"ahler--Einstein polarizations of a product of irreducible compact Hermitian symmetric spaces. Near the invariant metric, vanishing of the log-term forces the K\"ahler form to lie in its biholomorphic orbit. The argument combines the first variation of the Kempf distortion function with a central-positioning procedure, reducing rigidity to polynomiality of spectral multipliers of the Laplace operator and the Berezin transform. For integral K\"ahler--Einstein polarizations, a common-pole obstruction then rules out all infinitesimal deformations transverse to the automorphism orbit.
\end{abstract}
\maketitle

\section{Introduction}\label{sec-introduction}

Let $(L, h)$ be a positive Hermitian line bundle over a connected compact complex manifold $M$ of complex dimension $n$. By a theorem of Boutet de Monvel and Sj\"ostrand \cite{BoutetSjostrand}, the Szeg\H{o} kernel of the unit circle bundle $X$ of the dual bundle $(L^*, h^{-1})$ has the form $\Pi = u\, \varpi^{- n - 1} + v \log \varpi$, where $\varpi$ is (an almost-analytic extension of) a defining function of the unit disc bundle and $u$, $v$ are smooth functions on $X \times X$; we call $v$ the \emph{log-term} (see Section~\ref{sec-tyz} for the precise definition). Lu and Tian \cite{TianLu} showed that the vanishing of the log-term forces the coefficients $a_j$ of the Tian--Yau--Zelditch expansion of the Kempf distortion function (see Section~\ref{sec-tyz}) to vanish for $j > n$ \cite[Thm.~2.3]{TianLu}. They conjectured that, up to automorphisms, the Fubini--Study metric is the only metric on $\mathcal O(1) \to \mathbb{CP}^n$ with vanishing log-term \cite[p.~354]{TianLu}, proved this for $n = 1$ \cite[Thm.~1.1]{TianLu}, proved local rigidity in every dimension \cite[Thm.~1.2]{TianLu}, and showed that, generically, the log-term does not vanish \cite[Thm.~1.3]{TianLu}. Engli\v{s} and Zhang \cite{EnglisZhang} showed that the log-term vanishes for the invariant metric on every homogeneous line bundle over a compact irreducible Hermitian symmetric space \cite[Thms.~3.2--3.5]{EnglisZhang}, and Arezzo, Loi and Zuddas \cite{ArezzoLoiZuddas} showed more generally that the log-term vanishes for every regular quantization, i.e.\ whenever the Kempf distortion function is constant for all $m$ \cite[Thm.~1]{ArezzoLoiZuddas}. In particular, this holds over every compact homogeneous Hodge manifold \cite[Cor.~2]{ArezzoLoiZuddas}. 

The question addressed here is whether a rigidity statement comparable with that of Lu and Tian holds for compact Hermitian symmetric spaces. More precisely, Lu and Tian prove that a Hermitian metric on $\mathcal O(1) \to \mathbb{CP}^n$ that is sufficiently close to the standard metric in the $C^{2n + 4}$ norm, and whose associated Szeg\H{o} kernel has vanishing log-term, has curvature equal to the Fubini--Study form up to a holomorphic automorphism. Our Theorem~\ref{THM-quantizesymmetrygcd} extends this rigidity statement to compact symmetric products with integral K\"ahler--Einstein polarizations, which we shall now describe. Write $M = \prod_{k = 1}^s M_k$, where each $M_k$ is an irreducible Hermitian symmetric space of compact type of complex dimension $n_k$, so $n = \dim_{\C}(M) = \sum_{k = 1}^s n_k$. Let $L_k$ be the positive generator of the Picard group $\operatorname{Pic}(M_k)$, with invariant Hermitian metric $h_k$ and curvature $\omega_k = \Ric(h_k)$. Denote by $\gamma_k$ the genus of the noncompact dual of $M_k$, and put
\begin{equation*}
    d := \gcd(\gamma_1, \ldots, \gamma_s),\qquad \mathcal L^{\gcd} := \boxtimes_{k = 1}^s L_k^{\gamma_k/d}, \qquad h^{\gcd} := \boxtimes_{k = 1}^s h_k^{\gamma_k/d}
\end{equation*}
(see Subsection~\ref{subsec-products} for the precise notation). For a K\"ahler form $\omega \overset{loc}{=} \frac{\sqrt{-1}}{2\pi} \sum_{j, k} G_{j \bar k}\, dz_j \wedge d\overline{z_k}$ we denote by $\Ric_\omega \overset{loc}{:=} -\frac{\sqrt{-1}}{2\pi} \partial \bar\partial \log(\det(G_{j\bar k}))$ its Ricci form, so that $[\Ric_\omega] = c_1(M)$, and we say that $\omega$ is K\"ahler--Einstein with Einstein constant $\lambda \in \R$ if $\Ric_\omega = \lambda\, \omega$. Since $L_k^{\gamma_k} \cong K_{M_k}^{-1}$ \cite[p.~904]{EnglisZhang}, we have $K_M^{-1} \cong \boxtimes_k L_k^{\gamma_k} = (\mathcal L^{\gcd})^d$: thus $d$ is the Fano index of $M$, and $\mathcal L^{\gcd}$ is the unique primitive element of $\operatorname{Pic}(M) \cong \Z^s$ on the ray of $c_1(M)$. The curvature $\omega^{\gcd} := \Ric(h^{\gcd})$ is K\"ahler--Einstein with Einstein constant $d$; we call $(\mathcal L^{\gcd}, h^{\gcd})$ the \emph{K\"ahler--Einstein polarization} of $M$. For a Hermitian metric $h$ on $\mathcal L^{\gcd}$, we write $\logterm(h) \equiv 0$ when the log-term of the Szeg\H{o} kernel of the circle bundle of $((\mathcal L^{\gcd})^*, h^{-1})$ vanishes to infinite order along $X$ (see Section \ref{sec-tyz}). Our main result is the following.

\begin{theor}\label{THM-quantizesymmetrygcd}
Let $M = M_1 \times \cdots \times M_s$ be a product of irreducible Hermitian symmetric spaces of compact type $M_k$, $k = 1, \dots, s$, endowed with its K\"ahler--Einstein polarization $(\mathcal L^{\gcd}, h^{\gcd})$. There exists $\varepsilon = \varepsilon(M) \in \R^{> 0}$ with the following property: if $h = e^\psi h^{\gcd}$ is a Hermitian metric on $\mathcal{L}^{\gcd}$ such that $\logterm(h) \equiv 0$ and $\| \psi \|_{C^{2n + 4}} < \varepsilon$, then there exists a holomorphic automorphism $f$ of $M$ with $f^*(\Ric(h)) = \omega^{\gcd}$.
\end{theor}

The invariant metric itself has vanishing log-term (see Proposition \ref{PROP-homogeneouslogterm}), so Theorem~\ref{THM-quantizesymmetrygcd} characterizes it, up to automorphisms, among nearby metrics with this property. For $s = 1$, Theorem~\ref{THM-quantizesymmetrygcd} shows that every irreducible compact Hermitian symmetric space, polarized by the positive generator of its Picard group with its invariant metric, is locally rigid in the above sense. Moreover, our method actually shows that Theorem~\ref{THM-quantizesymmetrygcd} holds for every integral K\"ahler--Einstein polarization $\bigl((\mathcal L^{\gcd})^q, (h^{\gcd})^q\bigr)$, $q \in \N^*$, with $\omega^{\gcd}$ replaced by $q\, \omega^{\gcd}$ (see Remark \ref{RMK-integralKEpolarizationshold}). 

Its proof rests on two ingredients developed here. The \emph{analytic reduction} follows the scheme of \cite{TianLu}: the first variation of the Kempf distortion function (cf. \cite[Prop.~4.2]{TianLu}), the Taylor and elliptic estimates for the Tian--Yau--Zelditch coefficients (cf. \cite[Lemma~5.1, Cor.~3.1]{TianLu}), and a quantitative central positioning (cf. \cite[Lemma~7.2]{TianLu}). If local rigidity failed, these tools would produce a nonzero smooth function, orthogonal to the constants and to the first Laplace eigenspace, whose nonvanishing isotypical components have polynomial spectral multipliers. The \emph{spectral obstruction} shows that such polynomiality occurs only for the constant and first-eigenvalue modes, giving the required contradiction. This step is new: in \cite[\S6]{TianLu} the case of $\mathbb{CP}^n$ is handled by explicit computations with the monomial basis of $H^0(\mathbb{CP}^n, \mathcal O(m))$, which have no analogue in higher rank. The role of the polarization is essential in our method since, on an unrescaled product, factors with different genera can cancel poles in the rational Berezin eigenvalues (an explicit example on $\mathbb{CP}^1 \times \mathbb{CP}^3$ is given below).

The paper is organized as follows. Sections~\ref{sec-geometric}--\ref{sec-tyz} establish conventions, derive the variation formula, recall the definition of the log-term and the asymptotic results of Lu--Tian. Section~\ref{SEC-isotypicalcomponents} diagonalizes the variation on isotypical components, while Section~\ref{sec-central} develops the central-positioning and limiting arguments. Sections~\ref{sec-products}--\ref{sec-explicit} describe product spectra and prove the polynomiality criterion. Section~\ref{sec-rigidity} completes the rigidity proof.

\section{Geometric and analytical preliminaries}\label{sec-geometric}

We fix the conventions used throughout. The definitions in this section apply to any compact polarized K\"ahler manifold; homogeneity will be imposed only when needed.

Let $M$ be a connected compact complex manifold of complex dimension $n \geq 1$ with complex structure $J$, and let $\pi \colon L \to M$ be a positive holomorphic line bundle with a positive Hermitian metric $h_0$. The curvature form of $h_0$ and the associated K\"ahler metric are, for a local nonvanishing holomorphic section $\sigma$ of $L$,
\begin{equation*}
    \omega := \Ric(h_0) \overset{loc}{=} - \frac{\sqrt{-1}}{2\pi} \partial \bar\partial \log h_0(\sigma,\sigma),\qquad g(\cdot, \cdot) = \omega(\cdot, J \cdot).
\end{equation*}
Set $dV_0 := \omega^n/n!$. For $m \geq 1$, let $L^m$ be the $m$-th tensor power of $L$, with induced metric $h_0^m$, and put $d_m := \dim H^0(M,L^m)$. We use the $L^2$-Hermitian product
\begin{equation*}
    \langle s, s'\rangle := \int_M h_0^m(s, s')\, dV_0 \qquad(s, s' \in H^0(M, L^m)).
\end{equation*}
Write $(L^*, h_0^{-1})$ for the dual Hermitian line bundle and $\pi_{L^*} \colon L^* \to M$ for its projection. Unless otherwise indicated, $C^r$-norms, $r \in \N$, are computed with respect to $g$: namely, for $f \in C^r(M, \C)$ set
\begin{equation*}
    \|f\|_r := \sum_{j = 0}^r \sup_M| \nabla^j f |_g,
\end{equation*}
where $\nabla$ is the Levi-Civita connection of $g$. For later convenience, recall that for each $r \in \N$ there is $c_r \in \R^{\geq 1}$ such that $\|f_1 f_2\|_r \leq c_r \|f_1\|_r \|f_2\|_r$ for any $f_1,\, f_2 \in C^r(M, \C)$.

Given $\mathcal{X},\, \mathcal{Y}$ Fr\'echet spaces, $U \subset \mathcal{X}$ open, and a differentiable operator $\mathcal{A} \colon U \to \mathcal{Y}$, for $x \in U$ we denote the differential of $\mathcal{A}$ at $x$ by $D_x(\mathcal{A})$ and for $f \in \mathcal{X}$ we denote the directional derivative of $\mathcal{A}$ at $x$ in the direction of $f$ by $D_x(\mathcal{A})[f]$. This is the notation used, for instance, in Sections \ref{sec-kempf} and \ref{sec-tyz}.

Given $I \subseteq \N, | I | = \infty$, and $\{ x_m \}_{m \in I} \subset \C$, we write $x_m \in \C[m]$ (respectively $x_m \in \C(m)$) if there is $P \in \C[z]$ (respectively $R \in \C(z)$) such that $x_m = P(m)$ (respectively $x_m = R(m)$) for all $m \in I$; $P$ and $R$ are then unique. This is the notation used, for instance, in Sections \ref{sec-kempf} and \ref{SEC-isotypicalcomponents}.

\begin{rmk}\label{RMK-conventionLaplacian}\rm
We adopt the \emph{nonnegative} convention for the real Laplace--Beltrami operator of $g$ on functions, namely
\begin{align*}
    \Delta_g := -\operatorname{div}_g \grad_g .
\end{align*}
Let $(z_1, \dots, z_n)$ be local complex coordinates on $M$. With respect to this coordinate system, if $\omega = \frac{\sqrt{-1}}{2\pi} \sum_{j, k} G_{j \bar k}\, dz_j \wedge d\overline{z_k}$, $(G^{j \bar k})_{j, k} := \big( (G_{j \bar k})_{j, k}^T \big)^{-1}$, then for $\varphi \in C^\infty(M, \C)$
\begin{align*}
    \Delta_g \varphi = - 4\pi \sum_{j, k} G^{j \bar k} \frac{\partial^2 \varphi}{\partial z_j \partial \overline{z_k}} .
\end{align*}
Define the rescaled Laplacian $\Delta_\omega := \frac{1}{4\pi}\Delta_g$. Then for $\varphi \in C^\infty(M, \C)$
\begin{align*}
    \Delta_\omega \varphi = - \sum_{j, k} G^{j \bar k} \frac{\partial^2 \varphi}{\partial z_j \partial \overline{z_k}} ,
\end{align*}
and in particular it holds
\begin{equation}\label{EQ-Laplacianvolumeform}
    \frac{1}{(n - 1)!} \frac{\sqrt{-1}}{2\pi} \partial \bar\partial \varphi \wedge \omega^{n - 1} = - (\Delta_\omega \varphi)\, dV_0 .
\end{equation}
Indeed, at $p \in M$, choose normal holomorphic coordinates, in which
\begin{equation*}
    \omega_p = \frac{\sqrt{-1}}{2\pi} \sum_j dz_j \wedge d\bar z_j.
\end{equation*}
Only the diagonal terms of $\partial \bar\partial\varphi$ contribute to the wedge product with $\omega_p^{n - 1}$, giving
\begin{align*}
    \frac{\sqrt{-1}}{2\pi(n - 1)!}(\partial\bar\partial\varphi)_p \wedge \omega_p^{n - 1}
    & = \Big(\sum_j\frac{\partial^2\varphi}{\partial z_j \partial \bar z_j}(p)\Big)(dV_0)_p = -(\Delta_\omega \varphi)(p)(dV_0)_p ,
\end{align*}
as claimed. This choice of sign for $\Delta_\omega$ agrees with that of \cite{Gaudbook} (see, in particular, \cite[(1.11.16), (1.15.6), (3.2.4), (3.2.10)]{Gaudbook}), the rescaling by $\frac{1}{4\pi}$ being consistent with our normalization of the curvature of a line bundle (cf. \cite[(1.7.2)]{Gaudbook}). However, \cite{TianLu} adopts the \emph{nonpositive} convention, which corresponds to $-\Delta_\omega$ (see \cite[(3.8)]{TianLu}). This convention is visible, for instance, in the expression of the Laplacian of $\mathbb{CP}^n$ in \cite[p.~378]{TianLu} and in the relation $\Delta\xi = -(n + 1) \xi$ in the proof of \cite[p.~385, Thm.~1.2]{TianLu}.
\end{rmk}

\begin{rmk}\label{RMK-conventionsymmetricmetric}\rm
Let $M$ be an irreducible Hermitian symmetric space of compact type. The Picard group $\operatorname{Pic}(M)$ of $M$ is isomorphic to $\Z$ (see \cite[Sec.~14.2 and p.~520]{BorelHirze}): in our work, $L$ will be the positive generator of $\operatorname{Pic}(M)$ and $h_0$ is chosen so that $\omega$ is the unique invariant representative of $c_1(L)$, which is K{\"a}hler--Einstein with constant the genus of the noncompact dual of $M$ (see \cite[pp.~520--521]{BorelHirze}). This choice agrees with that of \cite{Zhang} (see \cite[p.~371 and pp.~377--378]{Zhang}).
\end{rmk}

\subsection{The Kempf distortion function}\label{subsec-kempf}

Assume $m$ is sufficiently large so that $L^m$ is very ample. Let $s_1, \dots, s_{d_m} \in H^0(M, L^m)$ be a $\langle \cdot, \cdot \rangle$-orthonormal basis. The \emph{Kempf distortion function} is $K_m(h_0) \in C^\infty(M, \R^{> 0})$ given by
\begin{equation}\label{DEF-Kempfdistortionfunction}
    K_m(h_0)(p) := \sum_{j = 1}^{d_m} h_0^m(s_j(p), s_j(p))\ ,\quad p \in M .
\end{equation}
This definition is independent of the choice of orthonormal basis, and very ampleness ensures strict positivity at every point. 
%More generally, for any positive Hermitian metric $h$ on $L$, we define $K_m(h)$ by the same formula, using $h^m$ and an orthonormal basis for the $L^2$-product determined by $h^m$ and the volume form $\Ric(h)^n/n!$.

\subsection{The Berezin transform}

We follow the approach of \cite{Schlichenmaier}, but we choose the opposite convention for Hermitian products, namely linearity in their first argument: accordingly, the corresponding identities have been adapted. Let $\Pi_m \colon L^2(M, L^m) \to H^0(M, L^m)$ be the orthogonal projection with respect to $\langle \cdot, \cdot \rangle$. For $\phi \in C^{\infty}(M, \C)$, we define the \emph{Toeplitz operator}
\begin{equation}\label{DEF-Toeplitzoperator}
    T^{(m)}_{\phi} \colon H^0(M, L^m) \to H^0(M, L^m) : s \mapsto \Pi_m(\phi \cdot s)
\end{equation}
If $L^m$ is very ample, the \emph{Berezin transform} $B_m(\phi) \in C^\infty(M, \C)$ of $\phi$ is defined as the function (see \cite[(6.5) and (7.1)]{Schlichenmaier})
\begin{equation}\label{DEF-Berezintransform}
    B_m(\phi)(p) = \frac{\langle T^{(m)}_{\phi} (e_{\alpha}^{(m)}), e_{\alpha}^{(m)} \rangle}{\langle e_{\alpha}^{(m)}, e_{\alpha}^{(m)} \rangle}\ ,\quad p \in M
\end{equation}
where, for every $\alpha \in \pi_{L^*}^{-1}(p) \smallsetminus \{ 0 \}$, $e_{\alpha}^{(m)} \in H^0(M, L^m)$ is uniquely determined by $\langle s, e_{\alpha}^{(m)} \rangle = \alpha^m(s(p))$ for every $s \in H^0(M, L^m)$, and we call it the \emph{coherent state vector} of $\alpha$. The quotient \eqref{DEF-Berezintransform} does not depend on the choice of $\alpha$ in the fiber over $p$, since replacing $\alpha$ by $c\alpha$, $c\in\C^*$, multiplies $e^{(m)}_\alpha$ by $\bar c^m$. The operator $B_m \colon C^\infty(M, \C) \to C^\infty(M, \C)$ is complex-linear.

The Berezin transform can be suitably expressed in terms of an orthonormal basis.

\begin{prop}\label{PROP-Berezin}
Let $s_1, \dots, s_{d_m} \in H^0(M, L^m)$ be a $\langle \cdot, \cdot \rangle$-orthonormal basis. Then for any $\phi \in C^{\infty}(M, \C)$ 
\begin{equation}\label{EQ-berezinpropdim}
    B_m(\phi) = \frac{1}{K_m(h_0)} \sum_{k, l = 1}^{d_m}\Big( \int_M \phi\ h_0^m(s_l, s_k)\, dV_0 \Big) h_0^m \left( s_k, s_l \right)
\end{equation}
\end{prop}
\begin{proof}
Let $p \in M$. Writing $e_\alpha^{(m)} = \sum_k c_k s_k$, $c_k \in \C$, orthonormality gives
\begin{equation*}
    B_m(\phi)(p) = \frac{\sum_{k, l} c_k \overline{c_l} \langle \phi s_k, s_l \rangle}{\sum_i |c_i|^2}.
\end{equation*}
The coherent-state identity \cite[Props.~6.5--6.6]{Schlichenmaier}, with the same convention, reads
\begin{equation}\label{EQ-factschlichen}
    h_0^m(s, s')(p) = \frac{\overline{\langle e_\alpha^{(m)},s\rangle} \langle e_\alpha^{(m)}, s'\rangle}{\langle e_\alpha^{(m)}, e_\alpha^{(m)} \rangle} K_m(h_0)(p).
\end{equation}
Since $\langle e_\alpha^{(m)}, s_j\rangle = c_j$, this identifies $\overline{c_l} c_k/ \sum_i| c_i|^2$ with
$h_0^m(s_l, s_k)(p)/K_m(h_0)(p)$. Substituting and using $\langle \phi s_k, s_l \rangle = \int_M \phi h_0^m(s_k, s_l)\, dV_0$ proves \eqref{EQ-berezinpropdim} after exchanging indices.
\end{proof}

\begin{rmk}\label{RMK-BerezinZhang}\rm
When $M$ is an irreducible compact Hermitian symmetric space and $(L, h_0)$ is as in Remark~\ref{RMK-conventionsymmetricmetric}, $B_m$ coincides with Zhang's Berezin transform $B_\nu$ for $\nu = m$. Indeed, \eqref{DEF-Berezintransform} is the covariant symbol of the Toeplitz operator $T^{(m)}_\phi$, i.e.\ Zhang's $B = W^*W$ \cite[\S1, (1.3)--(1.5)]{Zhang}; Proposition~\ref{PROP-Berezin} is the kernel form \cite[(2.13)]{Zhang}; and Zhang's line bundle $L$ is the generator of $\operatorname{Pic}(M)$, its power $L^\gamma$ being the anticanonical bundle \cite[\S2.2]{Zhang}. This identification is what allows us to import Zhang's eigenvalue formula in Section~\ref{sec-explicit}.
\end{rmk}

\section{First variation of the Kempf distortion function}\label{sec-kempf}

Fix $m$ such that $L^m$ is very ample, and choose a $\langle \cdot, \cdot \rangle$-orthonormal basis $s_1, \dots, s_{d_m}$ of $H^0(M,L^m)$ as in Subsection~\ref{subsec-kempf}. To compute the first variation of $K_m$ defined in \eqref{DEF-Kempfdistortionfunction}, write nearby positive Hermitian metrics as $h = e^\psi h_0$, with $\psi$ in the open set
\begin{align*}
    \mathcal{H}^\infty := \{ \psi \in C^\infty(M, \R) : \Ric(e^\psi h_0) > 0 \}
\end{align*}
in the Fr\'echet space $C^\infty(M, \R)$. The distortion function then defines a \emph{Kempf distortion operator}
\begin{align*}
    \mathcal{K}_m \colon \mathcal{H}^\infty \to C^\infty(M, \R) : \psi \mapsto K_m(e^\psi h_0) .
\end{align*}
We first express its dependence on $\psi$. Let $s^\psi_1, \dots, s^\psi_{d_m} \in H^0(M, L^m)$ be an orthonormal basis for the $L^2$-product $\langle \cdot, \cdot \rangle_\psi$ defined by $e^\psi h_0$, set $s^\psi_i = \sum_k a_{ki} s_k$, $a_{ki} \in \C$. Then
\begin{equation}\label{EQ-Gramrelation}
    \delta_{ij} = \langle s^\psi_i, s^\psi_j \rangle_\psi = \sum_{k, l} a_{k i} \overline{a_{lj}}\, \langle s_k, s_l \rangle_\psi .
\end{equation}
Let $G_\psi := (\langle s_k, s_l \rangle_\psi)_{k, l} \in M(d_m, \C)$ be the associated Gram matrix, $A := (a_{k i})_{k, i} \in \GL(d_m, \C)$, $I_{d_m} \in \GL(d_m, \C)$ the identity matrix. Then \eqref{EQ-Gramrelation} reads $I_{d_m} = A^T G_\psi \overline{A}$, so $G_\psi \in \GL(d_m, \C)$ and $G_\psi^{-1} = \overline{A} A^T$. Thus we compute
\begin{equation}\label{EQ-Kempfvariable}
\begin{split}
    \mathcal{K}_m(\psi) &= \sum_j (e^\psi h_0)^m(s^\psi_j, s^\psi_j) = e^{m \psi} \sum_{j, k, l} a_{k j} \overline{a_{lj}}\, h_0^m(s_k, s_l) \\
    &= e^{m \psi} \sum_{k, l} (G_\psi)^{-1}_{l k} h_0^m(s_k, s_l) = e^{m \psi} \Trace(G_\psi^{-1} H_0) ,
\end{split}
\end{equation}
where $H_0 := (h_0^m(s_k, s_l))_{k, l} \in C^\infty(M, M(d_m, \C))$. Observe that $\Trace(H_0) = \mathcal{K}_m(0)$. To study the variation of $\mathcal{K}_m$ at the origin, fix $\varphi \in C^\infty(M, \R)$, pick $t_0 \in (0, 1)$ such that $t\varphi \in \mathcal{H}^\infty$ for $t \in (-t_0, t_0)$, set $G_t := G_{t \varphi}$. For $t\neq0$, the difference quotient is
\begin{equation}\label{EQ-differencequotient}
    \frac{\mathcal{K}_m(t\varphi) - \mathcal{K}_m(0)}{t} = \frac{e^{m t \varphi} - 1}{t} \Trace(G_t^{-1}H_0) + \Trace\Big(\frac{G_t^{-1} - I_{d_m}}{t} H_0\Big) .
\end{equation}
We estimate both terms in each $C^r$-norm. Compactness of $M$ gives $C_{r,\varphi}>0$ such that
\begin{equation}\label{EQ-exponentialcontrol}
    \Big\| \frac{e^{m t \varphi} - 1}{t} - m \varphi \Big\|_r = m^2 |t| \Big\| \varphi^2 \int_0^1 (1 - \tau) e^{m \tau t\varphi} d\tau \Big\|_r \leq C_{r, \varphi} |t| .
\end{equation}
For the remaining terms, we need to control $G_t$ up to second order. Recall that 
\begin{equation}\label{EQ-Gramentries}
    (G_t)_{k l} = \int_M e^{m t \varphi} h_0^m(s_k, s_l) \frac{(\omega - t \frac{\sqrt{-1}}{2\pi} \partial \bar\partial \varphi)^n}{n!} .
\end{equation}
Observe that by the binomial theorem
\begin{equation}\label{EQ-volumebinomial}
\begin{split}
    \frac{(\omega - t \frac{\sqrt{-1}}{2\pi} \partial \bar\partial \varphi)^n}{n!} &= \sum_{\ell = 0}^n \frac{t^\ell}{\ell! (n - \ell)!} \Big(- \frac{\sqrt{-1}}{2\pi} \partial \bar\partial \varphi\Big)^\ell \wedge \omega^{n - \ell} \\
    &= dV_0 - \frac{t}{(n - 1)!} \frac{\sqrt{-1}}{2\pi} \partial \bar\partial \varphi \wedge \omega^{n - 1} + t^2 R(\varphi, t)\, dV_0 ,
\end{split}
\end{equation}
where, for $n \geq 2$, $R(\varphi, t) \in C^\infty(M,\R)$ is determined by the identity of top-degree forms
\begin{equation*}
    \sum_{\ell = 2}^n \frac{t^{\ell - 2}}{\ell! (n - \ell)!}
    \left(-\frac{\sqrt{-1}}{2\pi} \partial \bar\partial \varphi\right)^\ell
    \wedge \omega^{n - \ell} = R(\varphi, t)\,dV_0 ,
\end{equation*}
while for $n = 1$, set $R(\varphi,t) = 0$. By compactness of $M$, $\sup_{t \in [-1, 1]} \| R(\varphi, t) \|_r < \infty$. Inserting \eqref{EQ-Laplacianvolumeform} into \eqref{EQ-volumebinomial}, we get
\begin{equation}\label{EQ-volume2ndorderexpansion}
    \frac{(\omega - t \frac{\sqrt{-1}}{2\pi} \partial \bar\partial \varphi)^n}{n!} = (1 + t \Delta_\omega \varphi + t^2 R(\varphi, t))\, dV_0 .
\end{equation}
Writing $e^{m t \varphi} = 1 + m t \varphi + t^2 S(\varphi, t)$, compactness of $M$ gives $\sup_{t \in [-1, 1]} \| S(\varphi, t) \|_r < \infty$, and grouping terms of the same order in \eqref{EQ-Gramentries} yields
\begin{equation}\label{EQ-Gramentries2}
\begin{split}
    (G_t)_{k l} &= \int_M h_0^m(s_k, s_l)\, dV_0 + t \int_M h_0^m(s_k, s_l) (m \varphi + \Delta_\omega \varphi)\, dV_0 \\
    &\hspace{5mm} + t^2 \int_M h_0^m(s_k, s_l) T(\varphi, t)\, dV_0
\end{split}
\end{equation}
where to ease the notation we have set
\begin{align*}
    T(\varphi, t) &:= S(\varphi, t) + R(\varphi, t) + m \varphi \Delta_\omega \varphi \\
    &\hspace{5mm} + t (m \varphi R(\varphi, t) + S(\varphi, t) \Delta_\omega \varphi) + t^2 S(\varphi, t) R(\varphi, t)
\end{align*}
Define the matrices $B(\varphi), E(\varphi, t) \in M(d_m, \C)$ by
\begin{align*}
    B(\varphi)_{kl} & := \int_M h_0^m(s_k, s_l)  (m \varphi + \Delta_\omega \varphi)\, dV_0 ,\\
    E(\varphi, t)_{kl} & := \int_M h_0^m(s_k, s_l) T(\varphi, t)\, dV_0 .
\end{align*}
Then \eqref{EQ-Gramentries2} reads
\begin{equation}\label{EQ-Gram2ndorder}
    G_t = I_{d_m} + t B(\varphi) + t^2 E(\varphi, t) .
\end{equation}
On one hand, by construction $\sup_{t \in [-1, 1]} \| T(\varphi, t) \|_0 < \infty$, so the entries of $E(\varphi, t)$ remain bounded as $t \to 0$ and by \eqref{EQ-Gram2ndorder} we get with respect to the Euclidean topology
\begin{align*}
    \lim_{t \to 0} G_t = I_{d_m} ,
\end{align*}
which also gives $\lim_{t \to 0} G_t^{-1} = I_{d_m}$ with respect to the Euclidean topology. On the other hand, \eqref{EQ-Gram2ndorder} implies
\begin{equation}
    \lim_{t \to 0} \frac{G_t^{-1} - I_{d_m}}{t} = \lim_{t \to 0} - G_t^{-1} (B(\varphi) + t E(\varphi, t)) = - B(\varphi)
\end{equation}
with respect to the Euclidean topology. Since the entries of $G_t$ do not depend on the point of $M$, we deduce
\begin{align*}
    &\lim_{t \to 0} \Trace(G_t^{-1}H_0) = \Trace(H_0) ,\\ 
    &\lim_{t \to 0} \Trace\Big(\frac{G_t^{-1} - I_{d_m}}{t} H_0\Big) = - \Trace(B(\varphi)H_0)
\end{align*}
with respect to the $C^r$-norm. Combine these limits with \eqref{EQ-exponentialcontrol} to get with respect to the Fréchet topology
\begin{equation}\label{EQ-directionalderivativeKempf}
    \lim_{t \to 0} \frac{\mathcal{K}_m(t\varphi) - \mathcal{K}_m(0)}{t} = m \varphi \Trace(H_0) - \Trace(B(\varphi)H_0) .
\end{equation}
Since $\varphi$ is arbitrary, $\mathcal{K}_m$ admits directional derivatives at $0$. We now study the differentiability of $\mathcal{K}_m$. Observe that if we set
\begin{align*}
    L(\varphi) &:= m \varphi \Trace(H_0) - \Trace(B(\varphi)H_0) \\
    &= m \mathcal{K}_m(0) \varphi - \sum_{k, l = 1}^{d_m} \Big( \int_M h_0^m(s_l, s_k) (m \varphi + \Delta_\omega \varphi)\, dV_0 \Big) h_0^m(s_k, s_l) ,
\end{align*}
then $L \colon C^\infty(M, \R) \to C^\infty(M, \R) : \varphi \mapsto L(\varphi)$ is linear. It is also continuous: using the inequalities $\| f_1 f_2 \|_r \leq c_r \| f_1 \|_r \| f_2 \|_r$ and $\| \Delta_\omega \psi \|_0 \leq C_\omega \| \psi \|_2$, $C_\omega \in \R^{> 0}$, we get for $r \in \N$ 
\begin{align*}
    \| L(\psi) &\|_r \leq m c_r \| \mathcal{K}_m(0) \|_r \| \psi \|_r + \Vol_\omega(M) \\
    &\times \sum_{k, l = 1}^{d_m} (m \| \psi \|_0 + \| \Delta_\omega \psi \|_0) \| h_0^m(s_l, s_k) \|_0 \| h_0^m(s_k, s_l) \|_r \\
    &\leq \left( m c_r \| \mathcal{K}_m(0) \|_r + \Vol_\omega(M) (m + C_\omega) \sum_{k, l = 1}^{d_m} \| h_0^m(s_k, s_l) \|_r^2 \right) \| \psi \|_{\max\{r, 2\}} .
\end{align*}
The coefficient in parentheses on the right-hand side is finite because $M$ is compact, thus $L$ is continuous. We then show that $\mathcal{K}_m$ is differentiable at $0$ and $D_0(\mathcal K_m) = L$, by proving the following quadratic remainder estimate: for $r \in \N$, there exist $\delta \in (0, 1)$ and $\Theta_r \in \R^{> 0}$ such that if $\psi \in \mathcal{H}^\infty$, $\| \psi \|_{\max\{r, 2\}} < \delta$, then
\begin{align*}
    \| \mathcal{K}_m(\psi) - \mathcal{K}_m(0) - L(\psi) \|_r \leq \Theta_r \| \psi \|_{\max\{r, 2\}}^2 .
\end{align*}
Inspired by the above computations, we look for a uniform version of the previous estimates. For $\psi \in \mathcal{H}^\infty$, write
\begin{equation}\label{EQ-1storderKempf}
\begin{split}
    \mathcal{K}_m(\psi) - \mathcal{K}_m(0) - L(\psi) &= e^{m \psi} \Trace(G_\psi^{-1} H_0) + \Trace(B(\psi)H_0) \\
    &\hspace{5mm} - (1 + m \psi) \Trace(H_0) .
\end{split}
\end{equation}
By compactness of $M$, it suffices to control up to second order the first two summands of the right-hand side. For the exponential term, observe that
\begin{equation}\label{EQ-uniformexponentialcontrol}
\begin{split}
    \| e^{m \psi} - 1 - m \psi\|_r &= m^2 \Big\| \psi^2 \int_0^1 (1 - \tau) e^{m \tau \psi} d\tau \Big\|_r \\
    &\leq m^2 c_r^2 \Big\| \int_0^1 (1 - \tau) e^{m \tau \psi} d\tau \Big\|_r \| \psi \|_r^2 .
\end{split}
\end{equation}
Using $e^{m \tau \psi} - 1 = \sum_{k = 1}^\infty \frac{(m \tau)^k}{k!} \psi^k$, $c_r \geq 1$ and $\tau \in [0, 1]$, we compute
\begin{align*}
    \| e^{m \tau \psi} \|_r \leq 1 + \sum_{k = 1}^\infty \frac{(m \tau)^k}{k!} c_r^{k - 1} \| \psi \|_r^k = 1 + \frac{e^{m c_r \tau \| \psi \|_r} - 1}{c_r} \leq e^{m c_r \| \psi \|_r} ,
\end{align*}
so applying the Banach integral inequality to \eqref{EQ-uniformexponentialcontrol} gives
\begin{equation}\label{EQ-uniformexponentialcontrol2}
\begin{split}
    \| e^{m \psi} - 1 - m \psi\|_r &\leq m^2 c_r^2 \int_0^1 (1 - \tau) \| e^{m \tau \psi} \|_r d\tau\ \| \psi \|_r^2 \\
    &\leq \frac{m^2 c_r^2 e^{m c_r \| \psi \|_r}}{2} \| \psi \|_r^2 .
\end{split}
\end{equation}
Hence, if $\| \psi \|_{\max\{r, 2\}} \leq 1$, then \eqref{EQ-uniformexponentialcontrol2} yields the uniform estimate
\begin{equation}\label{EQ-uniformexponentialcontrol3}
    \| e^{m \psi} - 1 - m \psi\|_r \leq \frac{m^2 c_r^2 e^{m c_r}}{2} \| \psi \|_r^2 .
\end{equation}
For the remaining terms, we need to control $G_\psi$. As previously done, introduce
\begin{align*}
    &1)\ \mathcal{R}(\psi)\, dV_0 := \sum_{\ell = 2}^n \frac{1}{\ell! (n - \ell)!} \Big(- \frac{\sqrt{-1}}{2\pi} \partial \bar\partial \psi\Big)^\ell \wedge \omega^{n - \ell}\quad (\mathcal{R} := 0\ \text{if}\ n = 1)\\
    &2)\ \mathcal{S}(\psi) := e^{m \psi} - 1 - m \psi ,\\
    &3)\ \mathcal{T}(\psi) := \mathcal{S}(\psi) + \mathcal{R}(\psi) + m \psi \Delta_\omega \psi + (m \psi \mathcal{R}(\psi) + \mathcal{S}(\psi) \Delta_\omega \psi) + \mathcal{S}(\psi) \mathcal{R}(\psi) ,\\
    &4)\ \mathcal{E}(\psi) := \Big(  \int_M h_0^m(s_k, s_l) \mathcal{T}(\psi)\, dV_0 \Big)_{k, l}
\end{align*}
so that $\frac{(\omega - \frac{\sqrt{-1}}{2\pi} \partial \bar\partial \psi)^n}{n!} = (1 + \Delta_\omega \psi + \mathcal{R}(\psi))\, dV_0$ and $G_\psi = I_{d_m} + B(\psi) + \mathcal{E}(\psi)$. Moreover, since every summand in $\mathcal{R}(\psi)\, dV_0$ is homogeneous of degree $l \geq 2$ in $\partial \bar\partial \psi$, there is $C \in \R^{> 0}$ such that if $\| \psi \|_{\max\{r, 2\}} \leq 1$ 
\begin{equation}\label{EQ-uniformremaindercontrol}
    \| \mathcal{R}(\psi) \|_0 \leq \sum_{l = 2}^{n} A(n, l, \omega) \| \psi \|_2^l \leq C \| \psi \|_2^2 ,
\end{equation}
and combining \eqref{EQ-uniformremaindercontrol} with \eqref{EQ-uniformexponentialcontrol3} gives $C'' \in \R^{> 0}$ such that
\begin{equation}\label{EQ-uniform2ndordercontrol}
    | \mathcal{E}(\psi) | \leq C' \| \mathcal{T}(\psi) \|_0 \leq C'' \| \psi \|_2^2 
\end{equation}
where $| \cdot |$ is the operator norm. Since there is $C''' \in \R^{> 0}$ such that $| B(\psi) | \leq C''' \| \psi \|_2$, we get $\widetilde{C} \in \R^{> 0}$ such that if $\| \psi \|_{\max\{r, 2\}} \leq 1$
\begin{equation}\label{EQ-uniformGramPerturbation}
    | B(\psi) + \mathcal{E}(\psi) | \leq \widetilde{C}\, \| \psi \|_2 .
\end{equation}
Choose $\delta \in (0, 1)$ such that $\widetilde{C} \delta \leq \frac{1}{2}$, and restrict the potential $\psi$ to the open $\delta$-ball around $0$ in the $C^{\max\{r, 2\}}$-norm. Then $| B(\psi) + \mathcal{E}(\psi) | < \frac{1}{2}$, so the Neumann series $\sum_{\ell = 0}^\infty (-1)^\ell (B(\psi) + \mathcal{E}(\psi))^\ell$ converges absolutely and therefore realizes
\begin{equation}\label{EQ-Neumannseriesinverse}
    G_\psi^{-1} = \sum_{\ell = 0}^\infty (-1)^\ell (B(\psi) + \mathcal{E}(\psi))^\ell .
\end{equation}
Set $\mathcal{F}(\psi) = -\mathcal{E}(\psi) + \sum_{\ell = 2}^\infty (-1)^\ell (B(\psi) + \mathcal{E}(\psi))^\ell$, so $G_\psi^{-1} = I_{d_m} - B(\psi) + \mathcal{F}(\psi)$. By the remainder estimate \eqref{EQ-uniform2ndordercontrol}, the perturbation bound \eqref{EQ-uniformGramPerturbation}, and submultiplicativity, we obtain
\begin{align*}
    | \mathcal{F}(\psi) | &\leq | \mathcal{E}(\psi) | + \sum_{\ell = 2}^\infty |B(\psi) + \mathcal{E}(\psi)|^\ell \\
    &= | \mathcal{E}(\psi) | + \frac{| B(\psi) + \mathcal{E}(\psi) |^2}{1 - |B(\psi) + \mathcal{E}(\psi)|} \leq C'' \| \psi \|_2^2 + 2 \widetilde{C}^2 \| \psi \|_2^2 ,
\end{align*}
so there is $C_F \in \R^{> 0}$ such that $| \mathcal{F}(\psi) | \leq C_F \| \psi \|_2^2$. Finally, since \eqref{EQ-1storderKempf} can be rewritten as
\begin{equation}\label{EQ-1storderKempf2}
\begin{split}
    \mathcal{K}_m(\psi) &- \mathcal{K}_m(0) - L(\psi) = (1 + m \psi + \mathcal{S}(\psi)) \Trace(\mathcal{F}(\psi) H_0) \\
    &\hspace{5mm} - (m \psi + \mathcal{S}(\psi)) \Trace(B(\psi) H_0) + \mathcal{S}(\psi) \Trace(H_0) ,
\end{split}  
\end{equation}
combining the quadratic bound for $\mathcal F(\psi)$ with the identity \eqref{EQ-1storderKempf2} and the exponential estimate \eqref{EQ-uniformexponentialcontrol3} we obtain the desired estimate. Summing up, for $\psi \in C^\infty(M, \R)$ we obtain
\begin{equation}\label{EQ-differentialat0Kempf}
\begin{split}
    D_0(\mathcal{K}_m)[\psi] &= m \mathcal{K}_m(0) \psi\, - \sum_{k, l = 1}^{d_m} \Big( \int_M h_0^m(s_l, s_k) (m \psi + \Delta_\omega \psi)\, dV_0 \Big) h_0^m(s_k, s_l) .
\end{split}
\end{equation}
Proposition~\ref{PROP-Berezin} gives the key operator identity
\begin{equation}\label{EQ-differentialat0Kempf2}
    D_0(\mathcal{K}_m)[\psi] = \mathcal{K}_m(0)(m \psi - B_m(m \psi + \Delta_\omega \psi)) .
\end{equation}
For a homogeneous simply connected $(M, \omega)$, the Kempf distortion function is constant; integration yields
\begin{align*}
    \mathcal{K}_m(0) = \frac{d_m}{\Vol_\omega(M)} =: C_{(m)} \in \R^{> 0} .
\end{align*}
Since $c_1(M) > 0$, Kodaira vanishing and Hirzebruch--Riemann--Roch imply $d_m = \chi(M, L^m) \in \C[m]$ for sufficiently large $m$, hence $C_{(m)} \in \C[m]$. The variation formula becomes
\begin{equation}\label{EQ-differentialat0Kempfhomogeneous}
    D_0(\mathcal{K}_m)[\psi] = C_{(m)}(m \psi - B_m(m \psi + \Delta_\omega \psi)) .
\end{equation}
Any Hermitian symmetric space of compact type admits a holomorphic embedding into a finite-dimensional projective space (see \cite[Thm.~2]{DiScalaLoi}), so $L^m$ is very ample for any $m \in \N^*$ and the formula applies for every $m \geq 1$. This spectral form of the variation will be used in Section~\ref{SEC-isotypicalcomponents}.

\section{The Szeg\H{o} log-term and Tian--Yau--Zelditch expansion}\label{sec-tyz}

We begin by defining the Szeg\H{o} log-term. Since $M$ is compact and $(L, h_0)$ is a positive polarization, the unit circle bundle $X := \{ x \in L^* : h_0^{-1}(x, x) = 1 \}$ of $(L^*, h_0^{-1})$ is a compact strictly pseudoconvex Cauchy--Riemann manifold. Writing $d\theta$ for the $1$-form dual to the generator of the $S^1$-action on $X$, $dV_X := (\pi_{L^*})^*(\frac{\omega^n}{n!}) \wedge d\theta$ is a volume form on $X$. The \emph{Szeg\H{o} kernel} $\Pi$ of $X$ is the kernel of the orthogonal projection of $L^2(X, dV_X)$ onto the closed subspace of boundary values of holomorphic functions on $D := \{ x \in L^* : h_0^{-1}(x, x) < 1 \}$. By \cite[Cor.~1.7]{BoutetSjostrand}, there are $u,\, v \in C^\infty(X \times X, \C)$ determined up to functions vanishing to infinite order along $X$ such that
\begin{align*}
    \Pi = \frac{u}{(-\sqrt{-1}\, \varpi)^{n + 1}} + v \log(-\sqrt{-1}\, \varpi) ,
\end{align*}
where $\varpi$ is an almost-analytic extension of $x \in L^* \mapsto h_0^{-1}(x, x) - 1 \in \R$, and with abuse of notation $(-\sqrt{-1}\, \varpi)^{- n - 1}$ denotes the distributional limit of the functions $(-\sqrt{-1}\, \varpi + t)^{- n - 1}$ as $t \to 0^+$. We say that the log-term vanishes identically, and we write $\logterm(h_0) \equiv 0$, if $v$ vanishes to infinite order along $X$, i.e. it admits an identically vanishing representative.

We now recall the asymptotic results needed to convert the geometric vanishing condition into a statement about the linearized distortion function. For a positive polarization, the Kempf distortion function has the Tian--Yau--Zelditch expansion
\begin{equation}\label{EQ-TYZexpansion}
    \mathcal{K}_m(0) \sim \sum_{i = 0}^\infty\, a_i\, m^{n - i}\ ,\quad \{ a_j \}_{j \in \N} \subset C^\infty(M, \R)
\end{equation}
(see \cite[Thm.~1]{Zelditch}). More precisely, for any $k \geq 0$ and $\ell \geq 1$ there exists $C_{k, \ell} \in \R^{> 0}$ such that
\begin{equation}\label{EQ-TYZestimate}
    \Big\| \mathcal{K}_m(0) - \sum_{i = 0}^{\ell - 1} a_i\, m^{n - i} \Big\|_k \leq C_{k, \ell}\, m^{n - \ell} .
\end{equation}
For each $j \in \N$, $a_j$ is a universal polynomial of Weyl weight $2j$ in the curvature of $\omega$ and its covariant derivatives (see \cite[Thm.~1.1]{coeffLu}), and $a_0 \equiv 1$. In our discussion we vary the Hermitian metric: here and below, we write $a_j(\cdot, h)$ for the Tian--Yau--Zelditch coefficients of a Hermitian metric $h$ on $L$. We use the following two results:
\begin{itemize}
    \item if $\logterm(h) \equiv 0$, then $a_k(\cdot, h) \equiv 0$ for $k \geq n + 1$ (see \cite[Thm.~2.3]{TianLu});
    \item for any $k \in \N^*$ and $l > 2k + 2$ there exists $\kappa \in \R^{> 0}$, depending on $l$ and the $C^{l - 4}$-bound of the Riemannian curvature of $\omega$, such that if  $\varphi \in C^\infty(M, \R)$ satisfies $\frac{1}{2} \omega - \frac{\sqrt{-1}}{2\pi} \partial \bar\partial \varphi > 0$ and $a_k(\cdot, e^\varphi h_0) \equiv 0$, then
    \begin{equation}\label{EQ-vanishingcoeffestimate}
        \| \varphi \|_l \leq \kappa \| \varphi \|_{2k + 2}
    \end{equation}
    (see \cite[Cor.~3.1]{TianLu})
\end{itemize}
We also recall with which regularity the coefficients depend on the metric. For $l \geq 0$, set
\begin{align*}
    \mathcal{H}^l := \{ \psi \in C^{2l + 2}(M, \R) : \Ric(e^\psi h_0) > 0 \}
\end{align*}
which is open in the $C^{2l + 2}$ topology. Since the $l$-th coefficient has Weyl weight $2l$, there is a well-defined operator
\begin{align*}
    A_l \colon \mathcal{H}^l \to C^0(M, \R) : \psi \mapsto a_l(\cdot, e^\psi h_0) .
\end{align*}
The universal curvature-polynomial description makes $A_l$ smooth in the indicated Banach topologies. We shall use the following Taylor estimates from \cite[Lemma~5.1]{TianLu}.

\begin{prop}\label{PROP-TaylorestimateTYZcoefficients}
Let $l \in \N^*$. There exists $C \in \R^{> 0}$, depending only on $l$ and on the $C^{2l - 2}$-bound of the curvature of $\omega$ such that for every $\varphi \in C^{2l + 2}(M, \R)$ satisfying $\frac{1}{2} \omega - \frac{\sqrt{-1}}{2\pi} \partial \bar\partial \varphi > 0$ and $\| \varphi \|_{2l + 2} \leq 1$, for every $t \in [0, 1]$
\begin{equation}\label{EQ-2ndorderestimateTYZcoefficients}
    \left\| A_l(t \varphi) - A_l(0) - t\, D_0(A_l)[\varphi] \right\|_0 \leq C t^2 .
\end{equation}
\end{prop}

\begin{comment}
\begin{prop}\label{PROP-TaylorestimateTYZcoefficients}
Let $l \in \N^*$. If $\varphi \in C^{2l + 2}(M, \R)$ satisfies $\frac{1}{2} \omega - \frac{\sqrt{-1}}{2\pi} \partial \bar\partial \varphi > 0$ and $\| \varphi \|_{2l + 2} \leq 1$, then:
\begin{enumerate}
    \item there exists $C \in \R^{> 0}$ such that for every $t \in [0, 1]$
    \begin{equation}\label{EQ-2ndorderestimateTYZcoefficients}
        \left\| A_l(t \varphi) - A_l(0) - t\, D_0(A_l)[\varphi] \right\|_0 \leq C t^2
    \end{equation}
    \item there exists $U \subset \mathcal{H}^l$ neighborhood of $0$ and $C' \in \R^{> 0}$  such that for every $\psi \in U$ and $t \in [0, 1]$
    \begin{equation}\label{EQ-1storderestimateTYZcoefficients}
        \left\| A_l(t \varphi + \psi) - A_l(\psi) \right\|_0 \leq C' \| \varphi \|_{2l + 2} t
    \end{equation}
\end{enumerate}
The constants depend on $l$ and on the $C^{2l - 2}$-bound of the curvature of the fixed metric $\omega$.
\end{prop}
\end{comment}

The differentiated expansion follows from \cite[Thms.~4.1--4.2]{TianLu}.

\begin{prop}\label{PROP-differentiateTYZexpansion}
Let $\varphi \in C^\infty(M, \R)$. Then 
\begin{equation}\label{EQ-differentiateTYZexpansion}
    D_0(\mathcal{K}_m)[\varphi] \sim \sum_{i = 0}^\infty D_0(A_i)[\varphi]\, m^{n - i}
\end{equation}
in the sense that for any $k \in \N$ there is $\vartheta_k \in \R^{> 0}$ depending on $k$, $M$ and $\varphi$ such that
\begin{equation}\label{EQ-differentiateTYZestimate}
    \Big\| D_0(\mathcal{K}_m)[\varphi] - \sum_{i = 0}^k D_0(A_i)[\varphi]\, m^{n - i} \Big\|_0 \leq \vartheta_k\, m^{n - k - 1} .
\end{equation}
\end{prop}
Lu and Tian also state the corresponding $C^\infty$-asymptotic expansion; only the $C^0$ estimate above is needed here.

For ease of comparison, and with the normalizations now fixed, we record the precise local rigidity theorem of Lu and Tian for the projective space \cite[Thm.~1.2]{TianLu}, introduced in Section~\ref{sec-introduction}. Here $h_0$ denotes the standard Hermitian metric on the hyperplane bundle $\mathcal{O}(1) \to \mathbb{CP}^n$ whose curvature is the Fubini--Study metric $\omega_{\FS}$.

\begin{theor}\label{THM-LuTianMAIN}
There exists $\varepsilon \in \R^{> 0}$ such that for any $h = e^\psi h_0$ Hermitian metric on $\mathcal{O}(1)$ with $\logterm(h) \equiv 0$, if $\| \psi \|_{2n + 4} < \varepsilon$, then there exists a holomorphic automorphism $f$ of $\mathbb{CP}^n$ with $f^*(\Ric(h)) = \omega_{\FS}$.
\end{theor}
Theorem~\ref{THM-quantizesymmetrygcd} replaces the projective space by a compact symmetric product with its K\"ahler--Einstein polarization.

\section{Isotypical decomposition and spectral multipliers}\label{SEC-isotypicalcomponents}

To diagonalize the variation formula, we recall the relevant Peter--Weyl decomposition. Let $G$ be a compact Lie group with normalized Haar measure $\mu$. The action $G \times L^2(G, d\mu) \to L^2(G, d\mu)$ given by
\begin{align*}
    (\gamma f)(g) := f(\gamma^{-1} g)\ ,\quad f \in L^2(G, d\mu),\, g \in G ,
\end{align*}
yields an infinite-dimensional representation of $G$, called the \emph{regular representation}. It is unitary by the $G$-invariance of $\mu$, so by the Peter--Weyl theorem in its matrix-coefficient form (see \cite[Thm.~1.12]{Knappbook} or \cite[Thm.~35.5]{Etingof}) we have the Hilbert sum decomposition
\begin{align*}
    L^2(G, d\mu) \cong \widehat\bigoplus_{\rho \in \Sigma} (V_{\rho} \otimes V_\rho^*) \cong \widehat\bigoplus_{\rho \in \Sigma} V_{\rho}^{\oplus \dim(V_{\rho})} ,
\end{align*}
where $\Sigma$ denotes the set of isomorphism classes of irreducible unitary finite-dimensional representations of $G$. Here and below, $\widehat\bigoplus$ denotes the completed orthogonal Hilbert direct sum. Explicitly, the matrix-coefficient map is
\begin{align*}
    v \otimes \ell \in V_{\rho} \otimes V_\rho^* \mapsto \Big[ g \in G \mapsto \ell\bigl(\rho(g^{-1})v\bigr) \in \C \Big] \in L^2(G, d\mu) ,
\end{align*}
and $\dim(V_\rho^*) = \dim(V_\rho)$. The images of the Peter–Weyl summands consist of smooth functions, as the matrix coefficients of finite-dimensional representations are smooth (see \cite[pp.~176--179]{Etingof}).

Now, suppose $M$ is a homogeneous space given by $M = G/K$, $K$ a closed subgroup of $G$. Since $\mu$ descends to a $G$-invariant measure $\eta$ on $M$, we can consider $L^2(M, d\eta)$ and the action $G \times L^2(M, d\eta) \to L^2(M, d\eta)$ given by
\begin{align*}
    (\gamma f)(x) := f(\gamma^{-1} x)\ ,\quad f \in L^2(M, d\eta),\, x \in M .
\end{align*}
Then we have another unitary representation of $G$ and identifying functions on $M$ with right $K$-invariant functions on $G$, one obtains (see \cite[Ch.~V, Thm.~3.5]{Helgasonbook} or \cite[Ex.~35.7(i)]{Etingof})
\begin{align*}
    L^2(M, d\eta) \cong \widehat\bigoplus_{\rho \in \Sigma} (V_{\rho} \otimes (V_\rho^K)^*) \cong \widehat\bigoplus_{\rho \in \Sigma} V_\rho^{\oplus \dim(V_\rho^K)}
\end{align*}
where $V_\rho^K := \{ v \in V_\rho : \rho(k)(v) = v\ \forall\, k \in K \}$. Set $(V_\rho^*)^K := \{ \ell \in V_\rho^* : \ell \circ \rho(k) = \ell\,\ \forall\, k \in K \}$. We identify $(V_\rho^K)^*$ with $(V_\rho^*)^K$ via restriction to $V_\rho^K$. Indeed, since $\rho$ is unitary, every $K$-invariant linear functional on $V_\rho$ vanishes on $(V_\rho^K)^\perp$, and conversely, every linear functional on $V_\rho^K$ extends uniquely to a $K$-invariant linear functional on $V_\rho$ by setting it equal to zero on $(V_\rho^K)^\perp$. The matrix-coefficient map is
\begin{align*}
    v \otimes \ell \in V_{\rho} \otimes (V_\rho^*)^K \mapsto \Big[ gK \in M \mapsto \ell \bigl( \rho(g^{-1} )v \bigr) \in \C \Big] \in L^2(M, d\eta) ,
\end{align*}
and $\dim((V_\rho^K)^*) = \dim(V_\rho^K)$; the images of the summands still consist of smooth functions, as they are realized by the right $K$-invariant matrix coefficients. Define the \emph{isotypical components} of $L^2(M, d\eta)$ to be (here $\operatorname{Im}$ denotes the image through the matrix-coefficient map)
\begin{align*}
    W_\mu := \operatorname{Im}\Big(V_\mu^{\oplus \dim(V_\mu^K)} \Big) \subset C^\infty(M, \C) ,
\end{align*}
so the sum reads $L^2(M, d\eta) = \widehat\bigoplus_{\mu \in \Sigma} W_\mu$. In particular, each $W_{\mu}$ is invariant under the action $f \mapsto \gamma f$ and is a representation of $G$, not necessarily irreducible.

We now apply this decomposition to the Laplacian and Berezin transform. In order to relate the two setups, we assume that $G$ acts on $(M, \omega)$ by holomorphic isometries: indeed in this case $dV_0$ is $G$-invariant, so there is $C \in \R^{> 0}$ such that $dV_0 = C\, d\eta$, by transitivity of the action. In particular, $L^2(M, d\eta) = L^2(M, dV_0) =: L^2(M)$.

We begin the discussion with $\Delta_\omega \colon C^\infty(M, \C) \to C^\infty(M, \C)$ on $(M, \omega)$. Observe that since $G$ acts by holomorphic isometries, the $G$-action on $L^2(M)$ restricts to $C^\infty(M, \C)$. 

\begin{prop}\label{PROP-laplaciancommutes}
If $G$ acts on $(M, \omega)$ by holomorphic isometries, then $\Delta_\omega$ is $G$-equivariant.
\end{prop}
\begin{proof}
For $\gamma \in G$, the isometry condition gives $\grad_g(f \circ \gamma^{-1}) = \gamma_*(\grad_g f)$. Likewise, $\operatorname{div}_g(\gamma_*Y) = (\operatorname{div}_gY) \circ \gamma^{-1}$. Hence $\Delta_g(f\circ\gamma^{-1}) = (\Delta_gf) \circ \gamma^{-1}$, and likewise for $\Delta_\omega = \frac1{4\pi} \Delta_g$.
\end{proof}

We next consider the Berezin transform.

\begin{prop}\label{PROP-Berezincommutes}
If $M$ is simply connected and $G$ acts on $(M, \omega)$ by holomorphic isometries, then $B_m$ is $G$-equivariant.
\end{prop}
\begin{proof}
For $\gamma \in G$, the bundles $L$ and $(\gamma^{-1})^*L$
have the same first Chern class; simple connectedness and the K\"ahler condition give $H^1(M, \mathcal O_M) = 0$, so $L$ and $(\gamma^{-1})^*L$ are isomorphic. Since their metrics have the same curvature, the isomorphism may be chosen unitary. It induces a unitary map $U_\gamma$ on $H^0(M, L^m)$ with
\begin{equation*}
    h_0^m(U_\gamma s,U_\gamma t) = h_0^m(s,t) \circ \gamma^{-1}.
\end{equation*}
Apply \eqref{EQ-berezinpropdim} to the orthonormal basis $\{U_\gamma s_j\}$ and change variables using the invariance of $dV_0$. Both numerator and denominator transform by composition with $\gamma^{-1}$, giving $B_m(\phi \circ \gamma^{-1}) = B_m(\phi) \circ \gamma^{-1}$.
\end{proof}

Suppose now that $M$ is a compact Hermitian symmetric space. Then for all $\rho \in \Sigma$ it holds $\dim(V_\rho^K) \leq 1$ (see \cite[Thm.~3.1]{DOW} or \cite[Prop.~3.1]{Zhang}), and in particular, for any $\rho \in \Sigma$ such that $\dim(V_\rho^K) \neq 0$ we have
\begin{align*}
    W_\rho \cong V_\rho .
\end{align*}
Hence the isotypical components are irreducible. Consider a $G$-equivariant linear operator $D \colon C^\infty(M, \C) \to C^\infty(M, \C)$. For $\mu \in \Sigma$, denote by $P_\mu \colon L^2(M) \to W_\mu$ the orthogonal projection. If $\rho \in \Sigma$, then for all $\mu \in \Sigma$
\begin{align*}
    P_\mu \circ D |_{W_\rho} \colon W_\rho \to W_\mu\ \text{is linear and $G$-equivariant} .
\end{align*}
Since the isotypical components are irreducible, we can apply Schur's Lemma (see \cite[Lemma~11.9]{Etingof}) twice: first we get $P_\mu \circ D |_{W_\rho} = 0$ when $\rho \neq \mu$, hence $D(W_\rho) \subseteq W_\rho$, and secondly $D |_{W_\rho} \colon W_\rho \to W_\rho$ is a scalar multiple of $\Id_{W_\rho}$. By Propositions \ref{PROP-laplaciancommutes} and \ref{PROP-Berezincommutes}, we can apply this fact to $\Delta_\omega$ and $B_m$: therefore for any $\mu \in \Sigma$ there are $\lambda_\mu,\, \beta_\mu^{(m)} \in \C$ such that for any $f \in W_\mu$
\begin{equation}\label{EQ-laplacianisotypical}
    \Delta_\omega (f) = \lambda_\mu f
\end{equation}
and
\begin{equation}\label{EQ-Berezinisotypical}
    B_m(f) = \beta_\mu^{(m)} f .
\end{equation}
%In particular, $\lambda_\mu \in \R ^{\geq 0}$ since $\Delta_\omega$ is self-adjoint and nonnegative. 
We still denote by $D_0(\mathcal{K}_m)$ the complex-linear extension of the differential at $0$ of the Kempf distortion operator. Combine \eqref{EQ-laplacianisotypical} and \eqref{EQ-Berezinisotypical} with \eqref{EQ-differentialat0Kempfhomogeneous} to get for $f \in W_\mu$
\begin{equation}\label{EQ-differentialat0Kempfsymmetricisotypicalcomponent}
    D_0(\mathcal{K}_m)[f] = C_{(m)} \big(m - \beta_\mu^{(m)} (m + \lambda_\mu)\big) f .
\end{equation}
Set $F_\mu^{(m)} := C_{(m)} \big(m - \beta_\mu^{(m)} (m + \lambda_\mu)\big) \in \C$. For $\xi \in C^\infty(M, \C)$, let $\xi = \sum_{\mu \in \Sigma} \xi_\mu$ be its Hilbert sum decomposition in isotypical components. Since each isotypical component belongs to an eigenspace of $\Delta_\omega$, the Sobolev argument in the proof of \cite[Thm.~3.1]{Kazdan} shows that this sum converges in the $C^\infty$-topology. By linearity and continuity of $D_0(\mathcal{K}_m)$ we then get the Hilbert sum decomposition
\begin{equation}\label{EQ-differentialat0Kempfsymmetricisotypicaldecomposition}
    D_0(\mathcal{K}_m)[\xi] = \sum_{\mu \in \Sigma}\, F_\mu^{(m)} \xi_\mu ,
\end{equation}
which in particular converges in the $C^\infty$-topology. Finally, for later convenience observe that since the isotypical components are $L^2$-orthogonal, for every $\mu \in \Sigma$ we have $\langle D_0(\mathcal{K}_m)[\xi], \xi_\mu \rangle_{L^2(M)} = F_\mu^{(m)} \langle \xi_\mu, \xi_\mu \rangle_{L^2(M)}$, which reads
\begin{equation}\label{EQ-differentialisotypicalint}
    \int_M D_0(\mathcal{K}_m)[\xi]\, \overline{\xi_\mu}\, dV_0 = F_\mu^{(m)} \| \xi_\mu \|_{L^2(M)}^2 .
\end{equation}
\begin{rmk}\label{RMK-isotypicaleigenvaluedifferentialat0Kempf}\rm
For any $\mu \in \Sigma$ with $\dim(V_\mu^K) \neq 0$: $\beta_\mu^{(m)}$ lies in the field $\C(m)$ of rational functions of $m$. For irreducible $M$ this follows from \cite[Thm.~A]{Zhang}, and for products it then follows from \eqref{EQ-eigenvaluesKEproduct} in Section \ref{sec-products}. This together with $C_{(m)} \in \C[m]$ yields $F_\mu^{(m)} \in \C(m)$.
\end{rmk}

\section{Central positioning and the limiting criterion}\label{sec-central}

The central-positioning procedure removes infinitesimal automorphism directions before taking a deformation limit. We recall the relation between these directions and Laplace eigenfunctions. For a compact complex manifold $(M, J)$, the automorphism group $\Aut(M, J)$ is a complex Lie group whose real Lie algebra is $\mathfrak h := \{X \in \Gamma^\infty(M, TM) : \mathcal L_X J = 0\}$ \cite[Ch.~III, Thm.~1.1]{KobTransf}. On a compact K\"ahler manifold $(M, J, g, \omega)$, every $X \in \mathfrak h$ admits the Hodge decomposition \cite[Lemma~2.1.1]{Gaudbook}
\begin{equation}\label{EQ-realholomvectordecompose}
    X = X_H + \grad_g f + J \grad_g u ,
\end{equation}
where $X_H$ is the Riemannian dual of a harmonic one-form and $f,\, u \in C^\infty(M, \R)$ are uniquely determined by requiring zero mean. We call $f$ the \emph{real potential} and $\mathcal F := f + \sqrt{-1} u$ the \emph{complex potential} of $X$. With $X^{(1, 0)} := \frac{1}{2}(X - \sqrt{-1} J X)$, the identities needed below are \cite[Lemma~2.1.2 and Rmk.~2.1.2]{Gaudbook}
\begin{equation}\label{EQ-Liederivativepotentials}
    \mathcal L_X\omega = 2\sqrt{-1}\partial \bar\partial f, \qquad \mathcal L_{X^{(1, 0)}}\omega = \sqrt{-1}\partial \bar\partial \mathcal F.
\end{equation}
The Lie algebra $\mathfrak k := \{X \in \mathfrak{h} : \mathcal L_X g = 0\}$ of Killing fields consists precisely of the divergence-free fields of $\mathfrak{h}$; equivalently, $f = 0$ in \eqref{EQ-realholomvectordecompose} \cite[Prop.~2.2.1]{Gaudbook}. A Killing field is \emph{Hamiltonian} if and only if $X_H = 0$, in which case $X = J \grad_g u$. If $c_1(M) > 0$, every Killing field is Hamiltonian \cite[Rmk.~4.12 and Ex.~4.16]{Szekelyhidibook}. We are now ready to recall the following fundamental result due to Matsushima (see \cite[Thm.~3.6.2]{Gaudbook}):

\begin{theor}\label{THM-Matsushimadecomposition} 
Let $M$ be a compact K{\"a}hler--Einstein manifold of real dimension $2n$ with $c_1(M) > 0$. We have the direct sum decomposition $\mathfrak{h} = \mathfrak{k} \oplus J \mathfrak{k}$. Moreover, if $\mathfrak{h} \neq \{ 0 \}$, then $\lambda_1 = \frac{s}{4\pi n}$ is the first nonzero eigenvalue of $\Delta_\omega$, where $s \in \R^{> 0}$ is the scalar curvature of $(M,\, g)$, and if $\Lambda_1$ denotes the corresponding eigenspace, then
\begin{equation}\label{EQ-isomorphismMatsushima}
    \Lambda_1 \to \mathfrak{k} : \phi \mapsto J \grad_g \phi
\end{equation}
is a real-linear isomorphism.
\end{theor}

In the statements concerning holomorphic vector fields and central positioning, $\Lambda_1$ denotes the \emph{real} first eigenspace of the Laplacian. In the complex $L^2$-isotypical decomposition, the same notation denotes its complexification; orthogonality of a real function to $\Lambda_1$ is unaffected by this convention. We shall use the following immediate consequence of Matsushima's theorem.

\begin{cor}\label{COR-realpotentialeigenfunction}
Let $M$ be a compact K{\"a}hler--Einstein manifold, $c_1(M) > 0$ and $\mathfrak{h} \neq \{ 0 \}$. The real potential of a real holomorphic vector field belongs to $\Lambda_1$.
\end{cor}
\begin{proof}
By Theorem \ref{THM-Matsushimadecomposition}, write $X = J\grad_g f_1 + J( J \grad_g f_2) = - \grad_g f_2 + J\grad_g f_1$ with $f_i \in \Lambda_1$. Then $-f_2$ is the real potential, since $f_2$ has zero mean because $\Delta_\omega f_2 = \lambda_1 f_2$ with $\lambda_1 > 0$.
\end{proof}

We first construct a controlled representative transverse to the automorphism orbit, and then extract the polynomiality condition required for rigidity. Let $r \in \N$. For $\phi \in \mathcal{H}^r$, set
\begin{align*}
    \omega_\phi := \omega - \frac{\sqrt{-1}}{2\pi} \partial \bar\partial \phi \in \Gamma^{2r}(M, \Lambda^{1, 1}_\R T^*M) .
\end{align*}
If $M$ is compact, $\omega$ is K{\"a}hler--Einstein and $c_1(M) > 0$, then $F^*(\omega)$ is cohomologous to $\omega$ for any $F \in \Aut(M, J)$. Hence, there exists a unique $\Phi(F, \phi) \in C^{2r + 2}(M, \R)$ such that
\begin{align*}
    F^*(\omega_\phi) = \omega - \frac{\sqrt{-1}}{2\pi} \partial \bar\partial \Phi(F, \phi)\ ,\quad \int_M \Phi(F, \phi)\, dV_0 = 0 .
\end{align*}
Indeed, by the $\partial \bar\partial$-lemma there exists $\eta_F \in C^\infty(M, \R)$ such that $F^*(\omega) = \omega - \frac{\sqrt{-1}}{2\pi} \partial \bar\partial \eta_F$, so
\begin{equation}\label{EQ-explicitnormalizedpotential}
    \Phi(F, \phi) := \eta_F + \phi \circ F - \frac{\int_M (\eta_F + \phi \circ F)\, dV_0}{\Vol_\omega(M)} \in C^{2r + 2}(M, \R)
\end{equation}
is independent of the choice of $\eta_F$ and satisfies the required properties (uniqueness is just the maximum principle for pluriharmonic functions). Observe that $\Phi(F, \phi)$ has at least the same regularity as $\phi$.

The following quantitative form of central positioning (cf. \cite[Def.~6.1]{centralposition}) extends the method of \cite[Lemma~7.2]{TianLu}. Let \lq\lq$\perp$\rq\rq\, denote the orthogonal complement with respect to the $L^2$-product on $L^2(M, dV_0)$ and $p_{\Lambda_1} \colon L^2(M, dV_0) \to \Lambda_1$ denote the orthogonal projection onto $\Lambda_1$.

\begin{theor}\label{THM-centralpositioning}
Let $(M, \omega)$ be a compact K{\"a}hler--Einstein manifold, $c_1(M) > 0$ and $\mathfrak{h} \neq \{ 0 \}$. For any $r \in \N$, there are $\mathcal{U}_r \subset \mathcal{H}^r$ neighborhood of $0$ and $C_r \in \R^{> 0}$ depending on $r$, $M$, $\omega$, such that if $\phi \in \mathcal{U}_r$ then there exists $F \in \Aut(M, J)$ with
\begin{enumerate}
    \item[$(*1)$] $\Phi(F, \phi) \in \Lambda_1^\perp$ ;
    \item[$(*2)$] $\| \Phi(F, \phi) \|_{2r + 2} \leq C_r\, \| \phi \|_{2r + 2}$ .
\end{enumerate}
\end{theor}
\begin{proof}
For $X \in J\mathfrak{k}$ write $F_X = \exp(X)$, and let $f_X \in \Lambda_1$ be its real potential (Corollary~\ref{COR-realpotentialeigenfunction}). Integrating
\eqref{EQ-Liederivativepotentials} along $F_{tX}$ we can choose in \eqref{EQ-explicitnormalizedpotential} the potential $\eta_{F_X} = -4 \pi \int_0^1 f_X \circ F_{tX}\, dt$. Then by \eqref{EQ-explicitnormalizedpotential}, the condition $\Phi(F_X, \phi) \perp \Lambda_1$ is equivalent to
\begin{equation}\label{EQ-centralpositioning1}
    p_{\Lambda_1}\left(-4 \pi \int_0^1 f_X \circ F_{tX}\, dt + \phi \circ F_X \right) = 0 .
\end{equation}
Set
\begin{equation*}
    \mathcal B(X) = -4\pi p_{\Lambda_1} \int_0^1 f_X \circ F_{tX}\, dt\ , \qquad S_\phi(X) = p_{\Lambda_1}(\phi \circ F_X) .
\end{equation*}
Smooth dependence of the flows gives $\mathcal B(0) = 0$ and $D_0\mathcal B[X] = -4\pi f_X$, and Corollary \ref{COR-realpotentialeigenfunction}, Theorem \ref{THM-Matsushimadecomposition} make $X \in J\mathfrak{k} \mapsto f_X \in \Lambda_1$ a well-defined linear isomorphism. Write $\mathcal B(X) = D_0\mathcal B[X] + R(X)$, where $R(0) = 0, D_0R = 0$. Hence \eqref{EQ-centralpositioning1} is the fixed-point equation $T_\phi(X) = X$ for
\begin{equation*}
    T_\phi(X) := -(D_0 \mathcal B)^{-1}(R(X) + S_\phi(X)) .
\end{equation*}
Let $\Jac$ denote the Jacobian determinant. For an orthonormal basis $\{\theta_j\}$ of $\Lambda_1$, change of variables yields
\begin{equation*}
    \langle S_\phi(X), \theta_j \rangle_{L^2} = \int_M\phi\,(\theta_j \circ F_{-X}) \Jac(F_{-X})\, dV_0.
\end{equation*}
Fix a Euclidean norm $| \cdot |$ on $J\mathfrak{k}$. The factors multiplying $\phi$ have uniformly bounded $X$-derivatives on a small closed ball, so $\|S_\phi(B) - S_\phi(A)\|_{L^2} \leq C \|\phi\|_0 |B - A|$. As $D_0R = 0$, on a ball of radius $\delta$ we also have $\|R(B) - R(A)\|_{L^2} \leq C' \delta|B - A|$. Consequently,
\begin{equation*}
    |T_\phi(B) - T_\phi(A)| \leq (C_1\delta + C_2 \|\phi\|_0) |B - A|, \qquad |T_\phi(0)| \leq C_3 \|\phi\|_0 .
\end{equation*}
Choose $\delta \in \R^{> 0}$, then $\eta \in \R^{> 0}$, so that $C_1 \delta + C_2 \eta \leq \frac{1}{2}$ and $C_3 \eta \leq \frac{\delta}{2}$. For $\|\phi\|_0 \leq \eta$, $T_\phi$ is a contraction of the closed $\delta$-ball into itself. Its fixed point $Z$ satisfies
$|Z| \leq 2C_3 \|\phi\|_0$ and $\Phi(F_Z, \phi) \in\Lambda_1^\perp$.

Finite dimensionality and the smoothness of the flows on the corresponding compact parameter set give $\| f_Z \|_{2r + 2} \leq C_r'|Z|$ and for every $w \in C^{2r + 2}(M, \R)$
\begin{equation}\label{EQ-centralpositioning13}
    \|w \circ F_{tX} \|_{2r + 2} \leq C_r'' \| w \|_{2r + 2} \qquad (|X| \leq \delta,\ 0 \leq t \leq 1) .
\end{equation}
Using the normalized-potential formula \eqref{EQ-explicitnormalizedpotential} and the expression for $\eta_{F_Z}$ given above, then subtracting the mean, we obtain
\begin{equation}\label{EQ-centralpositioning15}
    \|\Phi(F_Z, \phi)\|_{2r + 2} \leq 2 C_r''(4\pi C_r'|Z| + \|\phi\|_{2r + 2}) \leq C_r \| \phi \|_{2r + 2}.
\end{equation}
The required neighborhood is $\mathcal U_r := \{\phi \in \mathcal H^r : \|\phi\|_{2r + 2} < \eta\}$.
\end{proof}

Central positioning and the estimates recalled in Section~\ref{sec-tyz} now give the analytic reduction used in the main theorem. In what follows, for $r, k \in \N$, we say that $\omega_1, \omega_2 \in \Gamma^r(M, \Lambda^k T^*M)$ are related if there is $F \in \Aut(M, J)$ such that $F^*(\omega_2) = \omega_1$. In this case, for the sake of brevity we just say $\omega_1, \omega_2$ are $F$-related. Also recall that for a Hermitian symmetric space of compact type: $\mathfrak{k} \neq \{ 0 \}$, so $\mathfrak{h} \neq \{ 0 \}$.

\begin{theor}\label{THM-polynomorthoglimit}
Let $(M, \omega = \Ric(h_0))$ be a Hermitian symmetric space of compact type with $\logterm(h_0) \equiv 0$ and $\omega$ K{\"a}hler--Einstein. If there exists $\{ \phi_i \}_{i \in \N} \subset C^{\infty}(M, \R)$ such that
\begin{enumerate}
    \item[(i)] $\| \phi_i \|_{2n + 4} \to 0$ as $i \to \infty$
    \item[(ii)] $\logterm(e^{\phi_i} h_0) \equiv 0$ 
    \item[(iii)] $\Ric(e^{\phi_i} h_0)$ is not $F$-related to $\omega$ for any $F \in \Aut(M, J)$
\end{enumerate}
then there exists $\xi \in (C^{\infty}(M, \R) \smallsetminus \{ 0 \}) \cap \Lambda_1^\perp$ such that $\int_M \xi\ dV_0 = 0$ and $F^{(m)}_{\mu} \in \C[m]$ for every $\mu \in \Sigma$ such that $\xi_{\mu} \neq 0$.
\end{theor}
\begin{proof}
Apply Theorem \ref{THM-centralpositioning} to obtain $u_i = \Phi(F_i, \phi_i) \in \Lambda_1^\perp$ with zero mean and $\| u_i \|_{2n + 4} \leq C_{n + 1} \| \phi_i \|_{2n + 4} \to 0$. By (iii), $u_i \neq 0$; biholomorphic invariance of the log-term implies $A_l(u_i) = 0$ for every $l \geq n + 1$. Put $\varepsilon_i = \| u_i \|_{2n + 4}$ and $\xi_i = u_i/\varepsilon_i$. The positivity of the nearby metrics and \eqref{EQ-vanishingcoeffestimate} give, for every $l > n + 1$, a uniform bound
$\| \xi_i \|_{2l + 2} \leq \kappa_l$. We also set $\kappa_{n + 1} := 1$, so the bound holds for $l \geq n + 1$. Compact embeddings and a diagonal argument yield a smooth limit
\begin{equation}\label{EQ-orderconvergence}
    \xi_i \longrightarrow \xi \quad\text{in}\ C^\infty(M, \R) \quad\text{along a subsequence}.
\end{equation}
Thus $\|\xi\|_{2n + 4} = 1$, while orthogonality and zero mean
pass to the limit.

For each $l \geq n + 1$, enlarge $\kappa_l$ if necessary so that $\frac{1}{2} \omega - \frac{\sqrt{-1}}{2\pi} \partial \bar\partial (\xi_i/\kappa_l) > 0$. Applying \eqref{EQ-2ndorderestimateTYZcoefficients} with $\varphi = \xi_i/\kappa_l$, $t = \kappa_l \varepsilon_i$ (which lies in $[0, 1]$ for large $i$) and using $A_l(u_i) = 0 = A_l(0)$ gives $D_0(A_l)[\xi_i] = O_l(\varepsilon_i)$. By continuity of $D_0(A_l)$ and \eqref{EQ-orderconvergence},
\begin{equation}\label{EQ-vanishingdifferentialcoefficientslimit}
    D_0(A_l)[\xi] = 0 \qquad(l \geq n + 1) .
\end{equation}
Proposition \ref{PROP-differentiateTYZexpansion} gives, for every $k \geq n$,
\begin{equation*}
    \| D_0(\mathcal K_m)[\xi] - P(\cdot, m) \|_0 \leq \vartheta_k m^{n - k - 1},\qquad P(\cdot, m) = \sum_{j = 0}^n D_0(A_j)[\xi] m^{n - j} .
\end{equation*}
For $\xi_\mu \neq 0$, let $P_\mu(m) = \langle P(\cdot, m), \xi_\mu \rangle_{L^2}/ \|\xi_\mu\|_{L^2}^2 \in \C[m]$. By \eqref{EQ-differentialisotypicalint}, $F_\mu^{(m)} - P_\mu(m) = O_{k, \mu}(m^{n - k - 1})$ for all $k \geq n$. But this difference is rational in $m$ (Remark \ref{RMK-isotypicaleigenvaluedifferentialat0Kempf}); a rational function decaying faster than any inverse power must vanish. Hence $F_\mu^{(m)} = P_\mu(m) \in \C[m]$.
\end{proof}

\section{Compact Hermitian symmetric spaces and their products}\label{sec-products}

Theorem~\ref{THM-polynomorthoglimit} reduces the rigidity problem to identifying the isotypical components for which $F_\mu^{(m)}=C_{(m)}\bigl(m-\beta_\mu^{(m)}(m+\lambda_\mu)\bigr)$ is polynomial in $m$. We prepare that classification by recording the spectral conventions on irreducible compact Hermitian symmetric spaces and the way they behave under products.

\subsection{Spherical representations and their parameters}

In Section~\ref{SEC-isotypicalcomponents} we have seen the isotypical decomposition of a  compact Hermitian symmetric space $M= G/K$
\begin{align*}
    L^2(M) = \widehat{\bigoplus}_{\rho \in \Sigma} W_\rho\ ,\quad W_\rho \cong V_\rho^{\oplus \dim(V_\rho^K)}\ ,\quad \dim(V_\rho^K) \leq 1.
\end{align*}
where $\Sigma$ denotes the set of the irreducible representation classes and the isotypical component $W_\rho$ is not trivial if and only if $V_{\rho}$ is {\it $K$-spherical}, i.e. $V_\rho^K \neq \{ 0 \}$. 

Here we recall a characterization of the spherical representations in terms of a root system of $G$ when $M$ is irreducible.

Let us first recall that, by the highest weight theorem (see \cite[Thm.~4.28]{Knappbook}), there is a bijection between $\Sigma$ and the set $P^+(G)$ of {\it dominant weights} of $G$.
More precisely, if $(\rho, V)$ is an irreducible finite-dimensional representation of $G$, then $\rho$ admits a unique highest weight $\mu \in P^+(G)$ (with respect to a choice of positive roots of $G$), and conversely every dominant weight $\mu$ occurs as the highest weight of an irreducible finite-dimensional representation $V_{\mu}$ of $G$. Setting $P_{\operatorname{sph}}^+(G, K) := \left\{ \mu\in P^+(G) : V_\mu^K \neq \{ 0 \} \right\}$, the isotypical decomposition reads
\begin{equation}\label{EQ-sphericaldecomposition}
    L^2(M) = \widehat{\bigoplus}_{\mu \in P_{\operatorname{sph}}^+(G, K)} W_\mu\ ,\quad W_\mu \cong V_\mu .
\end{equation}

For irreducible Hermitian symmetric spaces of compact type, $P_{\operatorname{sph}}^+(G, K)$ can be described in terms of a root system of $G$ by the Cartan--Helgason theorem (see \cite[Ch.~V, Sec.~4, Thm.~4.1]{Helgasonbook}) or equivalently as follows (see \cite[Sec.~2.1]{Zhang}). Let $r := \rank(M)$. We say two roots $\alpha, \beta$ are \emph{strongly orthogonal} if $\alpha + \beta,\, \alpha - \beta$ are not roots. One can construct a maximal set of strongly orthogonal noncompact positive roots $\alpha_1 < \dots < \alpha_r$, called the \emph{Harish-Chandra strongly orthogonal roots}, then take their images $\beta_1, \dots, \beta_r$ through the Cayley transform: by \cite[Prop.~3.1]{Zhang}, $\mu \in P_{\operatorname{sph}}^+(G, K)$ if and only if there are integers $\mu_1 \geq \mu_2 \geq \cdots \geq \mu_r \geq 0$ such that
\begin{align*}
    \mu = \sum_{j = 1}^r \mu_{r + 1 - j} \beta_j.
\end{align*}
Thus $P_{\operatorname{sph}}^+(G, K)$ is parametrized by 
\begin{align*}
    \mathcal{P}_r := \bigl\{ (m_1, \dots, m_r) \in \N^r : m_1 \geq \dots \geq m_r \bigr\}.
\end{align*}
For $\mathbf{m} = (m_1, \dots, m_r) \in \mathcal{P}_r$, denote by $W_{\mathbf{m}}$ the irreducible isotypical component associated with the highest weight $\sum_{j = 1}^r m_{r + 1 - j} \beta_j$. Then \eqref{EQ-sphericaldecomposition} becomes
\begin{equation}\label{EQ-HarishChandraisotypical}
    L^2(M) = \widehat{\bigoplus}_{\mathbf{m} \in \mathcal{P}_r} W_{\mathbf{m}}.
\end{equation}

This parametrization is all that is needed in the spectral calculations below.

\subsection{Products and tensor decompositions}\label{subsec-products}

We have completed the irreducible case, so we turn to the reducible setting. Let $M = G/K$, $\widetilde M = \widetilde G/\widetilde K$ be irreducible Hermitian symmetric spaces of compact type, $n := \dim_\C(M)$, $\widetilde n := \dim_\C(\widetilde M)$, $r := \rank(M)$, $\widetilde r := \rank(\widetilde M)$, $L$, $\widetilde L$ positive generators of $\operatorname{Pic}(M)$, $\operatorname{Pic}(\widetilde M)$, $h$, $\widetilde h$ positive Hermitian metrics on $L$, $\widetilde L$ such that $\omega = \Ric(h)$, $\widetilde \omega = \Ric(\widetilde h)$ are the unique $G$-invariant and $\widetilde G$-invariant forms with $[ \omega ] = c_1(L)$, $[ \widetilde \omega ] = c_1(\widetilde L)$, $g$, $\widetilde g$ the associated \K metrics, $p_1 \colon M \times \widetilde M \to M$, $p_2 \colon M \times \widetilde M \to \widetilde M$ the canonical projections. Consider the \K product $(M \times \widetilde M, g \times \widetilde g)$. The K{\"a}hler form of $g \times \widetilde g$ is
\begin{align*}
    p_1^*\omega + p_2^*\widetilde \omega =: \omega \boxtimes \widetilde \omega.
\end{align*}
Let $L \boxtimes \widetilde L := p_1^* L \otimes p_2^* \widetilde L$. For simplicity, denote it by $\mathcal{L}$. Observe that 
\begin{align*}
    c_1(\mathcal{L}) = p_1^* c_1(L) + p_2^* c_1(\widetilde L) = [ \omega \boxtimes \widetilde \omega ]
\end{align*}
and $h \boxtimes \widetilde h := p_1^* h \otimes p_2^*\widetilde h$ is a Hermitian metric on $\mathcal{L}$ such that
\begin{align*}
    \Ric(h \boxtimes \widetilde h) = \omega \boxtimes \widetilde \omega.
\end{align*}
Since $p_1^*\omega^{n + 1} = 0$, $p_2^*\widetilde \omega^{\widetilde n + 1} = 0$, by the binomial formula the K{\"a}hler volume form of $\omega \boxtimes \widetilde \omega$ satisfies
\begin{align*}
    dV_{\omega \boxtimes \widetilde \omega} = \frac{(\omega \boxtimes \widetilde \omega)^{n + \widetilde n}}{(n + \widetilde n)!} = \frac{p_1^*\omega^n}{n!} \wedge \frac{p_2^*\widetilde \omega^{\widetilde n}}{\widetilde n!} =: dV_\omega \boxtimes dV_{\widetilde\omega}.
\end{align*}
Since the volume form on the product is the product of the volume forms of the factors, the classical Hilbert tensor product theorem \cite[Sec.~II.4, Thm.~II.10(a)]{ReedSimon}  gives the canonical unitary identification
\begin{equation}\label{EQ-L2tensorproduct}
    L^2(M \times \widetilde M) \overset{U^{-1}}{\cong} L^2(M) \widehat\otimes L^2(\widetilde M) ,
\end{equation}
where $\widehat\otimes$ denotes the completed Hilbert tensor product and
\begin{align*}
    U(f \otimes \widetilde f)(x, \widetilde x) = f(x) \widetilde f(\widetilde x) =:(f \boxtimes \widetilde f)(x, \widetilde x).
\end{align*}
The hypotheses apply because compact manifolds with smooth volume measures have separable $L^2$ spaces. In particular,
\begin{align*}
    \langle f \boxtimes \widetilde f, u \boxtimes \widetilde u \rangle = \langle f, u \rangle_{L^2(M)} \langle \widetilde f, \widetilde u\rangle_{L^2(\widetilde M)} .
\end{align*}
Applying \eqref{EQ-HarishChandraisotypical} to both factors and distributing the Hilbert tensor product over orthogonal Hilbert sums yields
\begin{equation}\label{EQ-HilbertL2decompositionproduct}
    L^2(M \times \widetilde M) = \widehat{\bigoplus}_{(\mu, \nu) \in \mathcal P_r \times \mathcal P_{\widetilde r}} (W_\mu \boxtimes \widetilde W_\nu) ,
\end{equation}
where $W_\mu \boxtimes \widetilde W_\nu := U(W_\mu \otimes \widetilde W_\nu)$. Orthogonality and completeness follow from the corresponding properties of the decompositions of the factors. Finally, we verify that \eqref{EQ-HilbertL2decompositionproduct} is actually the isotypical decomposition under the action of $G \times \widetilde G$, which is given for $(a, b) \in G \times \widetilde G,\, F \in L^2(M \times \widetilde M)$ by
\begin{align*}
    ((a, b) \cdot F)(x,\widetilde x) = F(a^{-1} x, b^{-1} \widetilde x)\ ,\quad (x, \widetilde x) \in M \times \widetilde M
\end{align*} 
Observe that for $f \in L^2(M)$, $\widetilde f \in L^2(\widetilde M)$ it holds $(a, b) \cdot (f \boxtimes \widetilde f) = (a \cdot f) \boxtimes (b \cdot \widetilde f)$, so each $W_\mu \boxtimes \widetilde W_\nu$ is $(G \times \widetilde G)$-invariant. If $\rho_\mu \colon G \to \GL(W_\mu)$, $\widetilde \rho_\nu \colon \widetilde G \to \GL(\widetilde W_\nu)$ are the irreducible representations, this defines the representation
\begin{align*}
    \rho_{\mu, \nu} \colon G \times \widetilde G \to \GL(W_\mu \boxtimes \widetilde W_\nu) : \rho_{\mu, \nu}(a, b) := U \circ (\rho_\mu(a) \otimes \widetilde \rho_\nu(b)) \circ U^{-1} .
\end{align*}
By the classification of irreducible representations of direct products of compact groups (see \cite[Ch.~I, Sec.~4.4 and Thm.~6]{Vinberg}), this representation is irreducible, and the representations corresponding to distinct pairs $(\mu, \nu)$ are inequivalent. Hence, \eqref{EQ-HilbertL2decompositionproduct} is the isotypical decomposition of $L^2(M \times \widetilde M)$. For simplicity, write $W_{(\mu, \nu)} := W_\mu \boxtimes \widetilde W_\nu \subset C^\infty(M \times \widetilde M, \C)$, so that \eqref{EQ-HilbertL2decompositionproduct} reads
\begin{equation}\label{EQ-isotypicaldecompositionproduct}
    L^2(M \times \widetilde M) = \widehat{\bigoplus}_{ (\mu, \nu) \in \mathcal{P}_r \times \mathcal{P}_{\widetilde r}} W_{(\mu, \nu)}
\end{equation}
The decomposition of the volume form gives another relevant identification. Let $m \in \N^*$. Clearly $(h \boxtimes \widetilde h)^m = h^m \boxtimes \widetilde h^m$. Since Hermitian symmetric spaces of compact type are projective, we can apply the holomorphic K{\"u}nneth formula in degree $0$ (see \cite[Lemma~33.29.1]{StacksProject}) and get a canonical isomorphism
\begin{align*}
    H^0(M, L^m) \otimes H^0(\widetilde M, \widetilde L^m) \overset{\cong}{\to} H^0(M \times \widetilde M, \mathcal{L}^m) : s \otimes \widetilde s \mapsto p_1^*s \otimes p_2^*\widetilde s =: s \boxtimes \widetilde s
\end{align*}
which is unitary, since for $s, s' \in H^0(M, L^m)$, $t, t' \in H^0(\widetilde M, \widetilde L^m)$ we compute
\begin{align*}
    \langle s \boxtimes t, s' \boxtimes t' \rangle &= \int_{M \times \widetilde M} (h^m \boxtimes \widetilde h^m)(s \boxtimes t, s' \boxtimes t')\, dV_{\omega \boxtimes \widetilde \omega} \\
    &= \left( \int_M h^m(s, s')\, dV_\omega \right) \left( \int_{\widetilde M} \widetilde h^m(t, t')\, dV_{\widetilde \omega} \right) = \langle s, s' \rangle \langle t, t' \rangle
\end{align*}

We now turn to the operators. Since the isotypical components are generated by finitely many elements of the form $f \boxtimes \widetilde f$, by linearity it suffices to understand how the operators act on such elements.

As for $\Delta_{\omega \boxtimes \widetilde \omega}$, for $\phi \in C^\infty(M, \C)$ and $\widetilde \phi \in C^\infty(\widetilde M, \C)$ it holds (see \cite[Sec.~1.11, Lemma~1.18]{Grigoryan})
\begin{equation}\label{EQ-Laplacianproduct}
    \Delta_{\omega \boxtimes \widetilde \omega}(\phi \boxtimes \widetilde \phi) = \Delta_\omega(\phi) \boxtimes \widetilde \phi + \phi \boxtimes \Delta_{\widetilde \omega}(\widetilde \phi) .
\end{equation}
The Berezin transform needs a short computation. Set $d_m := \dim H^0(M, L^m)$ and $\widetilde d_m := \dim H^0(\widetilde M, \widetilde L^m)$. Choose orthonormal bases $\{s_k\}_{k = 1}^{d_m}$ of $H^0(M, L^m)$ and $\{\widetilde s_i\}_{i = 1}^{\widetilde d_m}$ of $H^0(\widetilde M, \widetilde L^m)$. By the K{\"u}nneth isomorphism, the sections
\begin{align*}
    S_{k, i} := s_k \boxtimes \widetilde s_i\ ,\quad 1 \leq k \leq d_m\ ,\quad 1 \leq i \leq \widetilde d_m ,
\end{align*}
form an orthonormal basis of $H^0(M \times \widetilde M, \mathcal{L}^m)$. Let $B_m$, $\widetilde B_m$, $B^\times_m$ and $K_m$, $\widetilde K_m$, $K^\times_m$ be the Berezin transforms and Kempf distortion functions respectively associated with $(L^m, h^m)$, $(\widetilde L^m, \widetilde h^m)$ and $(\mathcal{L}^m, (h \boxtimes \widetilde h)^m)$. Observe that for $(x, \widetilde x) \in M \times \widetilde M$
\begin{align*}
    K^\times_m(x, \widetilde x) &= \sum_{k = 1}^{d_m} \sum_{i = 1}^{\widetilde d_m} (h \boxtimes \widetilde h)^m(S_{k, i}, S_{k, i})(x, \widetilde x) \\ 
    &= \sum_{k = 1}^{d_m} \sum_{i = 1}^{\widetilde d_m} h^m(s_k, s_k)(x)\, \widetilde h^m(\widetilde s_i, \widetilde s_i)(\widetilde x) = K_m(x) \, \widetilde K_m(\widetilde x) ,
\end{align*}
which means $K^\times_m = K_m \boxtimes \widetilde K_m$. Then, by \eqref{EQ-berezinpropdim} we compute for $\phi \in C^\infty(M, \C)$ and $\widetilde \phi \in C^\infty(\widetilde M, \C)$
\begin{align*}
    \sum_{k, l = 1}^{d_m}& \sum_{i, j = 1}^{\widetilde d_m} \Big( \int_{M \times \widetilde M} (\phi \boxtimes \widetilde \phi)\, (h \boxtimes \widetilde h)^m(S_{lj}, S_{ki})\, dV_{\omega \boxtimes \widetilde\omega} \Big) \frac{(h \boxtimes \widetilde h)^m(S_{ki}, S_{lj})}{K^\times_m} \\
    &= \sum_{k, l = 1}^{d_m} \sum_{i, j = 1}^{\widetilde d_m} \left( \int_M \phi\, h^m(s_l, s_k)\, dV_\omega \right) \left( \int_{\widetilde M} \widetilde \phi\, \widetilde h^m(\widetilde s_j, \widetilde s_i)\, dV_{\widetilde\omega} \right) \\
    &\hspace{5mm}\times \frac{h^m(s_k, s_l) \boxtimes \widetilde h^m(\widetilde s_i, \widetilde s_j)}{K_m \boxtimes \widetilde K_m} ,
\end{align*}
which becomes
\begin{equation}\label{EQ-Berezinproduct}
    B^\times_m(\phi \boxtimes \widetilde\phi) = B_m(\phi) \boxtimes \widetilde B_m(\widetilde \phi) .
\end{equation}
The product formulas also determine the eigenvalues on each isotypical component. For $(\mu,\nu)\in\mathcal P_r\times\mathcal P_{\widetilde r}$, denote the Laplace and Berezin eigenvalues on $W_\mu$ by $\lambda_\mu$ and $\beta_\mu^{(m)}$, and those on $\widetilde W_\nu$ by $\widetilde\lambda_\nu$ and $\widetilde\beta_\nu^{(m)}$. Use $\lambda_{\mu,\nu}$ and $\beta_{\mu,\nu}^{(m)}$ for the respective eigenvalues on $W_{(\mu, \nu)}$. By \eqref{EQ-Laplacianproduct} and \eqref{EQ-Berezinproduct} we respectively obtain
\begin{equation}\label{EQ-laplacianeigenvalueproduct}
    \lambda_{\mu, \nu} = \lambda_\mu + \widetilde \lambda_\nu
\end{equation}
and
\begin{equation}\label{EQ-Berezineigenvalueproduct}
    \beta^{(m)}_{\mu, \nu} = \beta^{(m)}_\mu\, \widetilde \beta^{(m)}_\nu .
\end{equation}
Recall that in our setup  $K_m$, $\widetilde K_m$, $K^\times_m$ are constants, which we denote by $C_{(m)}$, $\widetilde C_{(m)}$, $C^\times_{(m)} \in \R^{> 0}$. Then $K^\times_m = K_m \boxtimes \widetilde K_m$ simply reads
\begin{equation}\label{EQ-Kempfproduct}
    C^\times_{(m)} = C_{(m)}\, \widetilde C_{(m)} .
\end{equation}
It follows from \eqref{EQ-differentialat0Kempfhomogeneous} that $D_0(\mathcal{K}^\times_m)$ acts on $W_{(\mu, \nu)}$ by scalar multiplication by
\begin{align*}
    F_{(\mu, \nu)}^{(m)} := C_{(m)}\, \widetilde C_{(m)} \big(m - \beta^{(m)}_\mu \cdot \widetilde \beta_{\nu}^{(m)}(m + \lambda_\mu + \widetilde \lambda_{\nu})\big) \in \C .
\end{align*}

\subsection{The integral K\"ahler--Einstein product polarization}

The preceding product formulas extend to $s\geq1$ irreducible compact Hermitian symmetric factors $M_1, \ldots, M_s$. For each $k$ write:
\begin{itemize}
    \item $M_k = G_k/K_k$, $M_k^*$ the noncompact dual of $M_k$
    \item $n_k := \dim_\C(M_k)$, $r_k := \rank(M_k)$, $\gamma_k$ the genus of $M_k^*$
    \item $L_k$ positive generator of $\operatorname{Pic}(M_k)$, $h_k$ positive Hermitian metric on $L_k$ such that $\omega_k = \Ric(h_k)$ is the unique $G_k$-invariant form with $[ \omega_k ] = c_1(L_k)$, set $C^k_{(m)} =  \frac{\dim(H^0(M_k, (L_k)^m))}{\Vol_{\omega_k}(M_k)} \in \R^{> 0}$
    \item $L^2(M_k) = \widehat \bigoplus_{\mu \in \mathcal{P}_{r_k}} W_{k, \mu}$ the decomposition \eqref{EQ-HarishChandraisotypical}
    \item for $\mu \in \mathcal{P}_{r_k}$, $\lambda_{k, \mu}$, $\beta^{(m)}_{k, \mu} \in \C$ the eigenvalues of $\Delta_{\omega_k}$, $B^k_m$ associated to $W_{k, \mu}$
\end{itemize}
Set $M := M_1 \times \dots \times M_s$, $\omega := \omega_1 \boxtimes \dots \boxtimes \omega_s$, $\mathcal{P} := \mathcal{P}_{r_1} \times \dots \times \mathcal{P}_{r_s}$, $\mathcal{L} := L_1 \boxtimes \dots \boxtimes L_s$, and consider the \K product $(M, \omega)$. The isotypical decomposition of $L^2(M)$ is given by
\begin{equation}\label{EQ-isotypicaldecompositiongeneralproduct}
    L^2(M) = \widehat \bigoplus_{(\mu^1, \dots, \mu^s) \in \mathcal{P}} W_{(\mu^1, \dots, \mu^s)}
\end{equation}
where $W_{(\mu^1, \dots, \mu^s)} = W_{1, \mu^1} \boxtimes \dots \boxtimes W_{s, \mu^s} \cong W_{1, \mu^1} \otimes \dots \otimes W_{s, \mu^s}$, and for $\mu = (\mu^1, \dots, \mu^s) \in \mathcal{P}$
\begin{equation}\label{EQ-parametersgeneralproduct}
    \lambda_{(\mu^1, \dots, \mu^s)} = \sum_{q = 1}^s \lambda_{q, \mu^q}\ ,\quad \beta_{(\mu^1, \dots, \mu^s)}^{(m)} = \prod_{j = 1}^s \beta_{j, \mu^j}^{(m)}
\end{equation}
while the Kempf distortion function is given by $\prod_{j = 1}^s C^j_{(m)} \in \R^{> 0}$.

For the rigidity theorem we use a common Einstein normalization. Set
\begin{align*}
    d := \gcd(\gamma_1, \dots, \gamma_s) \in \N^*
\end{align*}
and set $\mathcal{L}^{\gcd} := L_1^{\frac{\gamma_1}{d}} \boxtimes \dots \boxtimes L_s^{\frac{\gamma_s}{d}}$, $h^{\gcd} := h_1^{\frac{\gamma_1}{d}} \boxtimes \dots \boxtimes h_s^{\frac{\gamma_s}{d}}$. Consider 
\begin{align*}
    \omega^{\gcd} := \Ric(h^{\gcd}) = \frac{\gamma_1}{d} \omega_1 \boxtimes \dots \boxtimes \frac{\gamma_s}{d} \omega_s .
\end{align*}
Since each $\frac{\gamma_k}{d} \omega_k$, $k \in \{ 1, \dots, s\}$, is K{\"a}hler--Einstein with Einstein constant $\frac{d}{\gamma_k} \gamma_k = d$, so is $\omega^{\gcd}$. Since rescaling the metric does not change the isotypical decomposition of each $L^2(M_k)$, $k \in \{ 1, \dots, s\}$, the isotypical decomposition \eqref{EQ-isotypicaldecompositiongeneralproduct} for $M$ is unchanged. Indeed, putting $t_k := \frac{\gamma_k}{d}$, the volume forms satisfy
\begin{align*}
    dV_{\omega^{\gcd}} = \Big( \prod_{k = 1}^s \big(\frac{\gamma_k}{d}\big)^{n_k} \Big)\, dV_\omega ,
\end{align*}
so $L^2(M, dV_{\omega^{\gcd}}) = L^2(M, dV_\omega)$, and the two $L^2$-inner products differ by a positive scalar. Therefore, the underlying representation of $G_1 \times \cdots \times G_s$ on functions, its irreducible subrepresentations and the orthogonality of distinct isotypical components are the same, and the decomposition \eqref{EQ-isotypicaldecompositiongeneralproduct} still holds with exactly the same components. Denote by $\Delta_{\omega^{\gcd}}$, $B^{\gcd}_m$, $C^{\gcd}_{(m)}$ the rescaled real Laplace--Beltrami operator, the Berezin transform and the Kempf distortion function associated to $(\mathcal{L}^{\gcd}, h^{\gcd})$, and for $\mu = (\mu^1, \dots, \mu^s) \in \mathcal{P}$ denote by $\lambda^{\gcd}_\mu \in \C$ and $\beta_\mu^{\gcd, (m)} \in \C$ the eigenvalues of $\Delta_{\omega^{\gcd}}$ and $B^{\gcd}_m$ associated to $W_\mu$. We then compute
\begin{equation}\label{EQ-eigenvaluesKEproduct}
    \lambda^{\gcd}_\mu = \sum_{k = 1}^s \frac{d}{\gamma_k} \lambda_{k, \mu^k}\ , \quad \beta_\mu^{\gcd, (m)} = \prod_{k = 1}^s \beta_{k, \mu^k}^{(\frac{\gamma_k}{d} m)}
\end{equation}
and since $\Vol_{\omega^{\gcd}}(M) = \Big( \prod_{k = 1}^s \big(\frac{\gamma_k}{d}\big)^{n_k} \Big) \prod_{k = 1}^s \Vol_{\omega_k}(M_k)$, we get
\begin{equation}\label{EQ-KempfKEproduct}
    C^{\gcd}_{(m)} = \prod_{k = 1}^s \big(\frac{\gamma_k}{d}\big)^{-n_k} C^k_{(\frac{\gamma_k}{d} m)} \in \C[m].
\end{equation}
Thus, for every $\mu = (\mu^1, \dots, \mu^s) \in \mathcal{P}$, $D_0(\mathcal{K}^{\gcd}_m)$ acts on $W_\mu$ by scalar multiplication by
\begin{equation}\label{EQ-eigenvalueKEproduct}
    \Big( \prod_{k = 1}^s \big(\frac{\gamma_k}{d}\big)^{-n_k} C^k_{(\frac{\gamma_k}{d} m)} \Big) \Big( m - \Big( \prod_{k = 1}^s \beta_{k, \mu^k}^{(\frac{\gamma_k}{d} m)} \Big) \Big( m + \sum_{k = 1}^s \frac{d}{\gamma_k} \lambda_{k, \mu^k} \Big) \Big) \in \C ,
\end{equation}
which we denote by $F_\mu^{\gcd, (m)}$. In particular, $F_\mu^{\gcd, (m)} \in \C(m)$, and $F_\mu^{\gcd, (m)} \in \C[m]$ if and only if
\begin{equation}\label{EQ-polynomialconditiongcdproduct}
    \Big( \prod_{k = 1}^s C^k_{(\frac{\gamma_k}{d} m)} \beta_{k, \mu^k}^{(\frac{\gamma_k}{d} m)} \Big) \Big(m + \sum_{k = 1}^s \frac{d}{\gamma_k} \lambda_{k, \mu^k}\Big) \in \C[m] .
\end{equation}
The next section determines exactly when \eqref{EQ-polynomialconditiongcdproduct} holds.

\begin{rmk}\label{RMK-generalrescaledquantities}\rm
The identities \eqref{EQ-eigenvaluesKEproduct} and \eqref{EQ-KempfKEproduct} do not depend on the specific choice of rescaling: for any $t_1, \dots, t_s \in \N^*$, if we equip $\mathcal L^{\mathbf t} := \boxtimes_{k = 1}^s L_k^{t_k}$ with the invariant Hermitian metric $h^{\mathbf t} := \boxtimes_{k = 1}^s h_k^{t_k}$, then the eigenvalues of $\Delta_{\Ric(h^{\mathbf t})}$ and $B^{\mathbf t}_m$ associated to $W_\mu$ are
\begin{align*}
    \lambda^{\Ric(h^{\mathbf t})}_\mu = \sum_{k = 1}^s \frac{1}{t_k} \lambda_{k, \mu^k}\ , \quad \beta_\mu^{\mathbf t, (m)} = \prod_{k = 1}^s \beta_{k, \mu^k}^{(t_k m)}
\end{align*}
and the Kempf distortion function is given by
\begin{align*}
    C^{\mathbf t}_{(m)} = \prod_{k = 1}^s t_k^{-n_k} C^k_{(t_k m)} \in \C[m].
\end{align*}
\end{rmk}

\section{Spectral formulas and polynomiality}\label{sec-explicit}

The spectral classification needed in the rigidity argument is the following.

\begin{theor}\label{THM-gcdproductispolynomial}
Let $M = M_1 \times \dots \times M_s$ be the product of irreducible Hermitian symmetric spaces of compact type $M_k$, $k = 1, \dots, s$. For $\mu \in \mathcal{P}$, \eqref{EQ-polynomialconditiongcdproduct} holds if and only if either $\mu^k = 0$ for every $k = 1, \dots, s$ or exactly one $\mu^k$ equals $(1, 0, \dots, 0)$ and all the other ones are zero. In such cases, $F_\mu^{\gcd, (m)} \equiv 0$.
\end{theor}

We prove this using explicit formulas for the Kempf distortion, Laplace spectrum and Berezin eigenvalues of irreducible factors. The Berezin eigenvalue formula is recalled from \cite{Zhang}; we then analyze its rational poles. Let $M = G/K$ be an irreducible Hermitian symmetric space of compact type. Recall $n = \dim_\C(M)$, let $r = \rank(M)$, let $\gamma$ be the genus of the noncompact dual $M^*$ of $M$, let $a,\, b$ be the numerical invariants of $M^*$, related to $\gamma,\, n,\, r$ by 
\begin{equation}\label{EQ-symmetryinvariants}
    \gamma = (r - 1) a + b + 2\ ,\quad n = r + \frac{r(r - 1)}{2} a + r b
\end{equation}
Let $\Gamma$ denote the Gamma function. For $c \in \R$, $\boldsymbol{v} = (v_1, \dots, v_r) \in \R^r$, set $c + \boldsymbol{v} := (c + v_1, \dots, c + v_r) \in \R^r$. Define the \emph{Gindikin Gamma function} as
\begin{equation}\label{EQ-GindikinGamma}
    \Gamma_{M^*} (c + \boldsymbol{v}) = \prod_{j = 1}^r \Gamma \Big(c + v_j - \frac{a}{2} (j - 1) \Big)
\end{equation}
When writing $\Gamma_{M^*}(c)$, we mean $\Gamma_{M^*}(c + \boldsymbol{0})$.

\subsection{Berezin-transform eigenvalues}

We first record the finite-level vanishing of the Berezin eigenvalue: for fixed $m\geq1$ and $\mu=(\mu_1,\ldots,\mu_r)\in\mathcal{P}_r$ with $\mu_1>m$,
\begin{align*}
    \beta^{(m)}_\mu = 0\ ,\ \text{that is}\ B_m |_{W_\mu} \equiv 0
\end{align*}
Indeed, let $\phi \in W_\mu$. According to \eqref{DEF-Berezintransform}, it suffices to show that $T^{(m)}_{\phi} (e_{\alpha}^{(m)}) \equiv 0$, i.e. $\phi e_{\alpha}^{(m)}$ projects to the null vector in $H^0(M, L^m)$. Thus, it suffices to show for every $t \in H^0(M, L^m)$
\begin{equation}\label{EQ-orthogcoherentstate}
    \langle \phi e_{\alpha}^{(m)}, t \rangle = 0 
\end{equation} 
Recall that $H^0(M, L^m) \otimes \overline{H^0(M, L^m)}$ can be identified with a subspace of $L^2(M)$ as follows: for any $s \otimes \bar t \in H^0(M, L^m) \otimes \overline{H^0(M, L^m)}$, one defines a linear operator $A_{(s, \bar t)} \colon H^0(M, L^m) \rightarrow H^0(M, L^m)$ by
\begin{align*}
    A_{(s, \bar t)}(u) := \langle u, t \rangle s\ ,\quad u \in H^0(M, L^m)
\end{align*}
and for any linear operator $A \colon H^0(M, L^m) \rightarrow H^0(M, L^m)$, we define $\sigma(A) \in C^\infty(M, \C)$ by
\begin{align*}
    \sigma(A)(p) = \frac{\langle A(e_{\alpha}^{(m)}), e_{\alpha}^{(m)} \rangle}{\langle e_{\alpha}^{(m)}, e_{\alpha}^{(m)} \rangle}\ ,\quad p \in M ,
\end{align*}
where $\alpha \in L^* \smallsetminus \{ 0 \}$ satisfies $\pi_{L^*}(\alpha) = p$ and $e_\alpha^{(m)} \in H^0(M, L^m)$ is the coherent state vector of $\alpha$. Combining these two steps, we can associate to $s \otimes \bar t \in H^0(M, L^m) \otimes \overline{H^0(M, L^m)}$ a function $\Theta_{(s, \bar t)} \in C^\infty(M, \C)$, which using \eqref{EQ-factschlichen} is given at $p \in M$ simply by
\begin{align*}
    \frac{\langle A_{(s,\bar t)}(e_\alpha^{(m)}), e_\alpha^{(m)} \rangle}{\langle e_\alpha^{(m)},e_\alpha^{(m)}\rangle} = \frac{\langle e_\alpha^{(m)}, t\rangle \langle s, e_\alpha^{(m)} \rangle}{\langle e_\alpha^{(m)}, e_\alpha^{(m)} \rangle} = \frac{h_0^m(s, t)(p)}{C_{(m)}} .
\end{align*}
With respect to this construction, \cite[Thm.~4.2]{Zhang} states that
\begin{equation}\label{EQ-Zhangdecomposition}   
    H^0(M, L^m) \otimes \overline{H^0(M, L^m)} \cong \bigoplus_{\substack{\mu \in \mathcal{P}_r\\ \mu_1 \leq m}} W_\mu .
\end{equation}
Since the isotypical decomposition is orthogonal, we get for $s, t \in H^0(M, L^m)$  
\begin{align*}
    0 &= \langle \phi, \Theta_{(t, \bar s)}\rangle_{L^2(M)}
      = \frac{1}{C_{(m)}} \int_M \phi\, \overline{h_0^m(t, s)}\, dV_0 \\
      &= \frac{1}{C_{(m)}}\int_M h_0^m(\phi s, t)\, dV_0
      = \frac{1}{C_{(m)}} \langle \phi s, t\rangle ,
\end{align*}
which proves \eqref{EQ-orthogcoherentstate}. For $m\geq\mu_1$, Zhang's formula \cite[Thm.~A]{Zhang} gives the eigenvalues explicitly. In what follows, for $\mu = (\mu_1, \dots, \mu_r) \in \mathcal{P}_r$ set $\mu^* = (\mu_r, \mu_{r - 1}, \dots, \mu_1)$.
\begin{theor}\label{THM-eigenvalueBerezin}
Let $M$ be an irreducible Hermitian symmetric space of compact type. For every $\mu = (\mu_1, \dots, \mu_r) \in \mathcal{P}_r$, if $m \geq \mu_1$ then
\begin{equation}\label{EQ-eigenvalueBerezinGindikin}
    \beta^{(m)}_\mu = \frac{\Gamma_{M^*}(m + \gamma - \frac{n}r) \Gamma_{M^*}(m + \gamma)}{ \Gamma_{M^*}(m + \gamma + \mu) \Gamma_{M^*}(m + \gamma - \frac{n}r - \mu^*)}
\end{equation} 
while if $1 \leq m < \mu_1$, one has $\beta_\mu^{(m)} = 0$.
\end{theor}

\begin{comment}
Theorem~\ref{THM-eigenvalueBerezin} is \cite[Thm.~3.6]{Zhang} (= Thm.~A, pp.~372 and 384), via Remark~\ref{RMK-BerezinZhang}, with the following dictionary: Zhang's $\nu$ is our $m$, Zhang's genus $p$ is our $\gamma$, and Zhang's $\mathbf m = (m_1\leq\dots\leq m_r)$ is our $\mu^*$, so that Zhang's $\mathbf m^*$ is our $\mu$. Zhang's formula holds for all $\nu \geq 0$: for $\nu < m_r$, the factor $\Gamma(\nu + 1 - m_r)$ in the denominator makes it vanish (with $1/\Gamma(-k) = 0$ for $k \in \N$), so that the case distinction in Theorem~\ref{THM-eigenvalueBerezin} is automatic.
\end{comment}

We observe that \eqref{EQ-eigenvalueBerezinGindikin} yields $\beta^{(m)}_\mu \neq 0$ if $m \geq \mu_1$. For $j \in \{ 1, \dots, r \}$ set $c_j := m + \gamma - \frac{a}{2}(j - 1) \in \frac{1}{2}\Z$. Since $1 \leq j \leq r$ and $a, b \geq 0$, $m > 0$, by \eqref{EQ-symmetryinvariants}
\begin{align*}
    c_j &= m + (r - 1) a - \frac{a}{2}(j - 1) + b + 2 =\\
    &= m + (r - j)a + \frac{a}{2}(j - 1) + b + 2 > 0 ,
\end{align*}
so $c_j \in \frac{1}{2} \Z^{> 0}$, and by \eqref{EQ-GindikinGamma}, \eqref{EQ-eigenvalueBerezinGindikin} reads
\begin{equation}\label{EQ-eigenvalueBerezinGamma}
\begin{split}
    \beta^{(m)}_\mu = \prod_{j = 1}^r \frac{\Gamma \big(c_j - \frac{n}r \big) \Gamma( c_j )}{ \Gamma( c_j + \mu_j ) \Gamma \big(c_j - \frac{n}r - \mu_{r - j + 1} \big)} .
\end{split}
\end{equation}
Moreover, by \eqref{EQ-symmetryinvariants},
\begin{align*}
    c_j - \frac{n}{r} = m + \frac{r - j}{2} a + 1,\qquad c_j - \frac{n}{r}-\mu_{r - j + 1} = m - \mu_{r - j + 1} + \frac{r - j}{2} a + 1 .
\end{align*}
All Gamma arguments in \eqref{EQ-eigenvalueBerezinGamma} are therefore positive when $m \geq \mu_1$; hence $\beta_\mu^{(m)} > 0$. Regrouping the Gamma factors in \eqref{EQ-eigenvalueBerezinGamma} and replacing $j$ by
$r - j + 1$ in the quotient $\Gamma(c_j - n/r)/\Gamma(c_j - n/r - \mu_{r - j + 1})$, we obtain
\begin{equation}\label{EQ-eigenvalueBerezinfinal}
    \beta^{(m)}_\mu = \prod_{j = 1}^r \frac{\Gamma(m + 1 + (j - 1)\frac{a}{2}) \Gamma(m + \gamma - \frac{a}{2}(j - 1))}{ \Gamma(m + 1 - \mu_j + \frac{a}{2}(j - 1))  \Gamma(m + \gamma + \mu_j - \frac{a}{2}(j - 1))} .
\end{equation}

\subsection{Laplacian eigenvalues}

With our nonnegative Laplacian convention, the spherical eigenvalues take the following form.
\begin{theor}\label{THM-eigenvaluesLaplacian}
Let $M = G/K$ be an irreducible Hermitian symmetric space of compact type. For every $\mu = (\mu_1, \dots, \mu_r) \in \mathcal{P}_r$ 
\begin{equation}\label{EQ-eigenvalueLaplacian}
    \lambda_\mu = \sum_{j = 1}^r \mu_j \big( (r - j) a + b + 1 + \mu_j \big) \in \N .
\end{equation}
\end{theor}
\begin{proof}
Let $\widetilde g$ be the invariant metric induced by the negative of the Killing form on the $(-1)$-eigenspace of the Cartan involution of the Lie algebra of $G$. It has Einstein constant $1/2$ \cite[Thm.~7.73]{Besse}; hence irreducibility and comparison with the Riemannian Ricci tensor $\operatorname{Ric}^{\mathrm{Riem}}(g) = 2\pi\gamma g$ give $g = (4\pi \gamma)^{-1} \widetilde g$ and $\Delta_\omega = \gamma \Delta_{\widetilde g}$. For $\boldsymbol\mu = \sum_j\mu_{r + 1 - j}\beta_j$, Freudenthal's formula \cite[Appendix]{Urakawa} reads
\begin{equation}\label{EQ-normalizedeigenvalue}
    \widetilde\lambda_\mu = \langle 2\rho + \boldsymbol\mu, \boldsymbol\mu \rangle .
\end{equation}
Moreover, \cite[Sec.~3.1]{SahiZhang} gives
\begin{equation}\label{EQ-sumpositiverootsnumerical}
    2\rho = \sum_{j = 1}^r((j - 1) a + b + 1)\beta_j + \delta , \qquad \langle \delta, \beta_j \rangle = 0 .
\end{equation}
Pairwise orthogonality of the $\beta_j$ therefore yields
\begin{equation}\label{EQ-eigenvalueCayleyroots}
    \lambda_\mu = \gamma \sum_{j = 1}^r \mu_j((r - j) a + b + 1 + \mu_j) \langle \beta_{r + 1 - j}, \beta_{r + 1 - j} \rangle .
\end{equation}
The Cayley transform preserves the Killing form; the Harish-Chandra strongly orthogonal roots are long \cite[Sec.~2.1]{Alldridge}, with squared norm $1/h^\vee$, where $h^\vee$ is the dual Coxeter number \cite[Lemma~2.3]{BiquardGarciaPradaRubio}. Since $h^\vee = \gamma$ \cite[Sec.~2.2]{BiquardGarciaPradaRubio}, \eqref{EQ-eigenvalueLaplacian} follows.
\end{proof}

\begin{cor}\label{COR-firsteigenspace}
Let $M$ be an irreducible Hermitian symmetric space of compact type. Set $e_1 := (1, 0, \dots, 0) \in \mathcal{P}_r$. Then $\lambda_{e_1} = \gamma$ and $W_{e_1} = \Lambda_1$.
\end{cor}
\begin{proof}
By \eqref{EQ-eigenvalueLaplacian} and \eqref{EQ-symmetryinvariants}, $\lambda_{e_1} = \gamma$. For every nonzero $\mu \in \mathcal P_r$, $\mu_1 \geq 1$, and
\begin{equation*}
 \lambda_\mu \geq \mu_1 \bigl((r - 1) a + b + 1 + \mu_1 \bigr) \geq \gamma .
\end{equation*}
Equality holds precisely when $\mu=e_1$.
Since the decomposition \eqref{EQ-HarishChandraisotypical} exhausts $L^2(M)$, there are no additional Laplace modes: the inequality above therefore identifies the entire first eigenspace with $W_{e_1}$.
\end{proof}

Theorem~\ref{THM-eigenvaluesLaplacian} and Corollary~\ref{COR-firsteigenspace} also identify the isotypical decomposition of the eigenspace $\Lambda^{\gcd}_1$ of the first nonzero eigenvalue of $\Delta_{\omega^{\gcd}}$. For $q \in \{1, \dots, s\}$, denote by $\Psi_q \in \mathcal{P}$ the $s$-tuple such that $\mu^q = e_1^q = (1, 0, \dots ,0)$, $\mu^j = 0$ for $j \neq q$. If $\mu = (\mu^1, \dots, \mu^s) \in \mathcal{P}$, for $k \in \{1, \dots, s\}$ such that $\mu^k \neq 0$ then $\lambda_{k, \mu^k} \geq \gamma_k$ and equality holds if and only if $\mu^k = e_1^k$. By \eqref{EQ-eigenvaluesKEproduct}, for $\mu \neq 0$ we have
\begin{align*}
    \lambda^{\gcd}_\mu = \sum_{k = 1}^s \frac{d}{\gamma_k} \lambda_{k, \mu^k} \geq d \cdot \#(\{ j : \mu^j \neq 0 \}) \geq d
\end{align*}
and equality holds if and only if $\#(\{ j : \mu^j \neq 0 \}) = 1$ and the corresponding $\mu^k$ is equal to $e_1^k$. In other words,
\begin{equation}\label{EQ-firsteigenvaluegcdproduct}
    \lambda^{\gcd}_\mu = d \quad \iff \quad \mu = \Psi_k \text{ for some } k \in \{1, \dots, s\}
\end{equation}
In particular, the first nonzero eigenvalue of $\Delta_{\omega^{\gcd}}$ is $d$, and the isotypical decomposition reads
\begin{equation}\label{EQ-firsteigenspacegcdproduct}
    \Lambda^{\gcd}_1 = \bigoplus_{k = 1}^s W_{\Psi_k} .
\end{equation}

\subsection{The Kempf distortion polynomial}

We now determine $C_{(m)}$. In what follows, for $x \in \C$, $k \in \N^*$ denote the rising Pochhammer symbol by $(x)_k := x (x + 1) \cdots (x + k - 1)$, and for $s,\, t \in \R^{\geq 0}$, $s > \frac{(r - 1)}{2}a$, denote the generalized Pochhammer symbol by $((s))_t := \frac{\Gamma_{M^*}(s + t)}{\Gamma_{M^*}(s)}$.

\begin{prop}\label{PROP-Kempfirreduciblesymmetric}
Let $M$ be an irreducible Hermitian symmetric space of compact type. For $m \in \N^*$
\begin{equation}\label{EQ-Kempfgeneral}
    C_{(m)} = \prod_{j = 1}^r \left( m + 1 + \frac{(j - 1)}{2} a \right)_{\gamma - 1 - (j - 1) a} .
\end{equation}
\end{prop}
\begin{proof}
Homogeneity gives $C_{(m)} = d_m/\Vol_\omega(M)$. The weighted Bergman realization and reproducing kernel \cite[Sec.~3 and Lemma~3.1]{EnglisZhang} yield
\begin{equation}\label{EQ-dimensionPochhammer}
    d_m = \frac{((m + \gamma - \frac{n}{r}))_{\frac{n}{r}}}{((\gamma - \frac{n}{r}))_{\frac{n}{r}}} .
\end{equation}
The Hua polynomial \cite[Thm.~2.3]{YinLuRoos} is 
\begin{equation}\label{EQ-Huapolynomial}
    \chi(s) = \prod_{j = 1}^r \left(s + 1 + \frac{j - 1}{2}a \right)_{\gamma - 1 - (j - 1) a} .
\end{equation}
Since $\gamma - \frac{n}{r} = 1 + \frac{r - 1}{2} a$, the definition of the generalized Pochhammer symbol gives
\begin{equation}\label{EQ-HuaHilbert}
    d_m = \frac{\chi(m)}{\chi(0)} .
\end{equation}
By \eqref{EQ-symmetryinvariants}, $\chi$ is monic of degree $\sum_j(1 + b + (r - j) a) = n$. Comparing its leading coefficient with $d_m = m^n \Vol_\omega(M) + O(m^{n - 1})$ from Hirzebruch--Riemann--Roch gives $\Vol_\omega(M) = 1/ \chi(0)$, hence $C_{(m)} = \chi(m)$.
\end{proof}

\begin{rmk}\label{RMK-zerosHua}\rm
By \eqref{EQ-Huapolynomial}, the zeros of $\chi$ are the numbers $- \bigl(1 + \frac{a}{2}(j - 1) + i \bigr)$ with $1 \leq j \leq r$ and $0 \leq i \leq \gamma - 2 - (j - 1)a$. Since $1 + \frac{a}{2}(j - 1) + \gamma - 2 - (j - 1) a = \gamma - 1 - \frac{a}{2}(j - 1)$, all of them lie in $[-(\gamma - 1), -1]$.
\end{rmk}

\subsection{The pole obstruction and the polynomiality criterion}

The obstruction to polynomiality is already visible in each factor of \eqref{EQ-eigenvalueBerezinfinal}. For $\mu_j = 0$ the corresponding factor of \eqref{EQ-eigenvalueBerezinfinal} equals $1$. For $\mu_j > 0$, the recurrence $\Gamma(z + 1) = z\Gamma(z)$ reduces that factor to
\begin{equation}\label{EQ-factoreigenvalueBerezinindex}
    \prod_{\ell = 0}^{\mu_j - 1} \frac{m + \frac{a}{2} (j - 1) - \ell}{m + \gamma - \frac{a}{2} (j - 1) + \ell} .
\end{equation}
We first establish the criterion for one irreducible factor.

\begin{prop}\label{PROP-productispolynomialirreducible}
Let $M$ be an irreducible Hermitian symmetric space of compact type. Then for any $\mu = (\mu_1, \dots, \mu_r) \in \mathcal{P}_r$ it holds
\begin{equation}\label{EQ-goalirreduciblecase}
    C_{(m)} \beta_\mu^{(m)} (m + \lambda_\mu) \in \C[m] \quad\iff\quad \mu \in \{ 0, e_1 \}
\end{equation}
and in both cases, $F_\mu^{(m)} \equiv 0$.
\end{prop}
\begin{proof}
For $\mu = 0$ or $\mu = e_1$, respectively, $(\beta_\mu^{(m)}, \lambda_\mu) = (1, 0)$ or
$(\beta_\mu^{(m)}, \lambda_\mu) = (m/(m + \gamma), \gamma)$. In either case $\beta_\mu^{(m)}(m + \lambda_\mu) = m$ and $F_\mu^{(m)} = 0$.
  
For $\mu \notin \{0, e_1\}$, let $\widehat\beta_\mu(z)$ be the rational extension of \eqref{EQ-factoreigenvalueBerezinindex}, $\widehat C(z) = \chi(z)$, and $R_\mu(z) = \widehat C(z) \widehat\beta_\mu(z)(z + \lambda_\mu)$. As $\mu_1 \geq 1$, $\widehat\beta_\mu$ has a denominator factor $z + \gamma$, while no numerator factor vanishes at $z = -\gamma$, since $\frac{a}{2}(j - 1) < \gamma$. Furthermore, $\widehat C(-\gamma) \neq 0$ by Remark \ref{RMK-zerosHua} and $\lambda_\mu > \gamma$ by Corollary \ref{COR-firsteigenspace}. Thus $R_\mu$ retains a pole at $-\gamma$, so
\begin{equation}\label{EQ-refinedgoalirreduciblecase}
    R_\mu \notin \C[z] .
\end{equation}
If $C_{(m)} \beta_\mu^{(m)} (m + \lambda_\mu)$ were a polynomial, it would coincide with $R_\mu$ on all integers $m \geq \mu_1$ and hence as a rational function, contradicting \eqref{EQ-refinedgoalirreduciblecase}.
\end{proof}

For unrescaled products with different genera, cancellations may invalidate this criterion. Consider $\mathbb{CP}^1 \times \mathbb{CP}^3$, so $\gamma_1 = 2$, $\gamma_2 = 4$. One has
\begin{align*}
    C^1_{(m)} = m + 1\ ,\quad C^2_{(m)} = (m + 1)(m + 2)(m + 3)
\end{align*}
Take $\mu = (\mu^1, \mu^2) = ((2), 0) \in \mathcal{P}$. By \eqref{EQ-eigenvalueBerezinGindikin} and \eqref{EQ-eigenvalueLaplacian}, we obtain
\begin{align*}
    &\beta_{1,\mu^1}^{(m)} = \frac{m(m - 1)}{(m + 2)(m + 3)} , &\lambda_{1,\mu^1} = 6 ,\\
    &\beta_{2,\mu^2}^{(m)} = 1 , &\lambda_{2,\mu^2} = 0 .
\end{align*}
Consequently, we compute
\begin{align*}
    &\Big( \prod_{j = 1}^2 C^j_{(m)} \beta_{j, \mu^j}^{(m)} \Big) \Big(m + \sum_{q = 1}^2 \lambda_{q, \mu^q}\Big) \\
    &\quad = (m + 1)^2(m + 2)(m + 3)
    \frac{m(m - 1)}{(m + 2)(m + 3)}(m + 6) \\
    &\quad = m(m - 1)(m + 1)^2(m + 6) \in \C[m] ,
\end{align*}
and $F_\mu^{(m)} = 12m(m + 1)^2 \in \C[m]$, although $\mu^1 \neq 0, e_1$. Thus the product polarization in Theorem~\ref{THM-gcdproductispolynomial} is not an inessential normalization: it prevents this cross-factor cancellation.

We can now prove the product criterion, Theorem \ref{THM-gcdproductispolynomial}, by locating the common pole.

\begin{proof}[Proof of Theorem~\ref{THM-gcdproductispolynomial}]
For $\mu = 0$ or $\mu = \Psi_k$, respectively, $(\beta_\mu^{\gcd,(m)}, \lambda_\mu^{\gcd}) = (1, 0)$ or $(\beta_\mu^{\gcd,(m)}, \lambda_\mu^{\gcd}) = (m/(m + d), d)$ by \eqref{EQ-eigenvaluesKEproduct} and \eqref{EQ-firsteigenvaluegcdproduct}. In either case $\beta_\mu^{\gcd,(m)}(m + \lambda^{\gcd}_\mu) = m$ and $F_\mu^{\gcd, (m)} = 0$.

Conversely, suppose \eqref{EQ-polynomialconditiongcdproduct} holds and $\mu \neq 0$. With the rational extensions notation used in Proposition \ref{PROP-productispolynomialirreducible}, set
\begin{equation*}
    R_\mu^{\gcd}(z) = (z + \lambda_\mu^{\gcd}) \prod_{j = 1}^s \widehat C^j \left( \frac{\gamma_j}{d} z \right) \widehat\beta_{j,\mu^j} \left( \frac{\gamma_j}{d} z \right) .
\end{equation*}
Since $R_\mu^{\gcd}(m) = \Big( \prod_{k = 1}^s C^k_{(\frac{\gamma_k}{d} m)} \beta_{k, \mu^k}^{(\frac{\gamma_k}{d} m)} \Big) \Big(m + \sum_{k = 1}^s \frac{d}{\gamma_k} \lambda_{k, \mu^k}\Big)$ for every $m \in \N^*$, the polynomiality assumption forces $R_\mu^{\gcd} \in \C[z]$. We now analyze the factors at $z = - d$:

\emph{(a) Kempf factors}. By Remark \ref{RMK-zerosHua}, each $\widehat C^k \left( \frac{\gamma_k}{d} z \right)$ vanishes only for $z \in [ -d \frac{\gamma_k - 1}{\gamma_k}, - \frac{d}{\gamma_k}] \subset (-d, 0)$, so $\widehat C^k \left( \frac{\gamma_k}{d} (-d) \right) \neq 0$.

\emph{(b) Berezin numerators.} Let $a^{(k)}, b^{(k)}$ denote the numerical invariants of $M_k^*$. By \eqref{EQ-factoreigenvalueBerezinindex}, the numerator of $\widehat\beta_{k, \mu^k}(\frac{\gamma_k}{d} z)$ is a product of factors $\frac{\gamma_k}{d} z + \frac{a^{(k)}}{2} (j - 1) - l$ with $l \geq 0$, which vanish at $z = \bigl(l - \frac{a^{(k)}}{2} (j - 1)\bigr)/\frac{\gamma_k}{d} \geq - \frac{a^{(k)}(r_k - 1)}{2\frac{\gamma_k}{d}} > - \frac{\gamma_k}{\frac{\gamma_k}{d}} = -d$, because $\frac{a^{(k)}}{2}(r_k - 1) < (r_k - 1) a^{(k)} + b^{(k)} + 2 = \gamma_k$. Hence no numerator factor vanishes at $z = -d$.

\emph{(c) Berezin denominators.} If $\mu^k \neq 0$, then $\mu^k_1 \geq 1$ and the denominator of $\widehat\beta_{k, \mu^k}(\frac{\gamma_k}{d} z)$ contains the factor $\frac{\gamma_k}{d} z + \gamma_k = \frac{\gamma_k}{d}(z + d)$, which vanishes at $z = -d$.

Since there is $k \in \{ 1, \dots, s \}$ with $\mu^k \neq 0$, (a)--(c) show that $\prod_{j = 1}^s \widehat C^j \left( \frac{\gamma_j}{d} z \right) \widehat\beta_{j,\mu^j} \left( \frac{\gamma_j}{d} z \right)$ has a pole at $z = -d$ of order at least $1$. Multiplication by $z + \lambda_\mu^{\gcd}$ lowers the order of this pole by at most one, and only if $\lambda_\mu^{\gcd} = d$. Since $R^{\gcd}_\mu$ is a polynomial, we must have $\lambda_\mu^{\gcd} = d$, hence $\mu = \Psi_k$ by \eqref{EQ-firsteigenvaluegcdproduct}.
\end{proof}

\begin{rmk}\label{RMK-generalrescalingpolynomiality}\rm
The same argument applies more generally, namely whenever $t_1, \dots, t_s \in \N^*$ satisfy $\frac{\gamma_1}{t_1} = \dots = \frac{\gamma_s}{t_s} = \rho \in \mathbb Q^{>0}$. Indeed, the common pole of the rescaled Berezin eigenvalues then occurs at $m = -\rho$, and the first nonzero eigenvalue of the corresponding K{\"a}hler--Einstein metric is $\rho$. Hence the conclusion of Theorem \ref{THM-gcdproductispolynomial} holds for every integral product polarization whose associated metric is K{\"a}hler--Einstein.
\end{rmk}

\section{Proof of the local rigidity theorem}\label{sec-rigidity}

We observe that the invariant reference polarization has vanishing log-term. This places it within the scope of the limiting criterion. Indeed, \cite[Cor.~2]{ArezzoLoiZuddas} applied to our setup, together with Remark \ref{RMK-generalrescaledquantities}, gives:

\begin{prop}\label{PROP-homogeneouslogterm}
Let $M = \prod_{k = 1}^s M_k$ be a product of irreducible Hermitian symmetric spaces of compact type. For any $t_1, \dots, t_s \in \N^*$, equip $\mathcal L^{\mathbf t} := \boxtimes_{k = 1}^s L_k^{t_k}$ with the invariant Hermitian metric $h^{\mathbf t} := \boxtimes_{k = 1}^s h_k^{t_k}$. Then the Kempf distortion function of $(\mathcal L^{\mathbf t}, h^{\mathbf t})$ is a polynomial in the quantization parameter $m$, and
\begin{equation*}\label{EQ-regularproductlogterm}
    \logterm(h^{\mathbf t}) \equiv 0 .
\end{equation*}
\end{prop}

The analytic reduction of Theorem~\ref{THM-polynomorthoglimit} and the pole obstruction of Theorem~\ref{THM-gcdproductispolynomial} now give the main result.

\begin{proof}[Proof of Theorem~\ref{THM-quantizesymmetrygcd}]
Otherwise, a sequence $\{\phi_i\}_i \subset C^\infty(M, \R)$ with $C^{2n + 4}$-norm converging to $0$ would define Hermitian metrics $\{ e^{\phi_i} h^{\gcd} \}_i$ with vanishing log-term without $\Ric(e^{\phi_i} h^{\gcd})$ being related to $\omega^{\gcd}$. This metric is K\"ahler--Einstein, and Proposition \ref{PROP-homogeneouslogterm} gives $\logterm(h^{\gcd}) \equiv 0$. Theorem \ref{THM-polynomorthoglimit} then provides $0 \neq \xi \in(\Lambda_1^{\gcd})^\perp$, of zero mean, with $F_\mu^{\gcd, (m)} \in \C[m]$ for every $\xi_\mu \neq 0$. By Theorem \ref{THM-gcdproductispolynomial} and \eqref{EQ-firsteigenspacegcdproduct}, 
\begin{equation}\label{EQ-specialgcdsubspaces}
    \xi \in W_0 \oplus \bigoplus_{k = 1}^s W_{\Psi_k} = W_0 \oplus \Lambda_1^{\gcd} .
\end{equation}
The zero-mean condition eliminates the constants $W_0$, while orthogonality eliminates $\Lambda_1^{\gcd}$, contradicting $\xi \neq 0$.
\end{proof}

\begin{rmk}\label{RMK-integralKEpolarizationshold}\rm
The proof of Theorem \ref{THM-quantizesymmetrygcd} applies verbatim to the integral K{\"a}hler--Einstein polarizations $\bigl( (\mathcal L^{\gcd})^q, (h^{\gcd})^q \bigr)$, $q \in \N^*$, with $\omega^{\gcd}$ replaced by $q\, \omega^{\gcd}$. Indeed:
\begin{itemize}
    \item by Remark \ref{RMK-generalrescaledquantities}, the relevant spectral data are those of $t_k = q \frac{\gamma_k}{d}$;
    \item by Remark \ref{RMK-generalrescalingpolynomiality} the conclusion of Theorem \ref{THM-gcdproductispolynomial} holds for them;
    \item by Proposition \ref{PROP-homogeneouslogterm} the invariant metric has vanishing log-term;
    \item the first eigenspace of $\Delta_{q \omega^{\gcd}}$ is again $\bigoplus_k W_{\Psi_k}$, with eigenvalue $\frac{d}{q}$.
\end{itemize}
This justifies the statement made after Theorem \ref{THM-quantizesymmetrygcd}.
\end{rmk}

\section*{Acknowledgements}
The authors  acknowledge the support of the Fondazione di Sardegna
through the project ProBiki, and of INdAM--GNSAGA
(Gruppo Nazionale per le Strutture Algebriche, Geometriche e loro
Applicazioni).

\section*{Declaration on the use of AI tools}
The ideas, results and proofs presented in this paper are entirely
due to the authors. Generative AI tools, namely Claude (Anthropic)
and ChatGPT (OpenAI), were used to assist in checking the formal
correctness of the arguments and in improving the presentation and
language of the manuscript. All outputs of these tools were
critically reviewed and verified by the authors, who take full
responsibility for the content of this paper.

\end{document}